\documentclass[12pt]{article}
\usepackage{amsmath}
\usepackage{mathdots}
\usepackage{amsfonts}
\usepackage{amssymb}
\usepackage{latexsym}
\usepackage{graphicx}
\usepackage{pgfplots}
\usepackage{todonotes}

\usepackage{color}
\usepackage{fullpage}
\usepackage{verbatim}
\usepackage[section]{algorithm}
\usepackage{algorithmic}

\def\eepsilon{\epsilon_*}
\def\qed{\hfill $\square $}

\def\O{\mathcal{O}}

\def\Pi{\prod}
\def\P{\mathcal{P}}

\def\xh{\widehat x}

\def\Xt{\widetilde X}

\def\Vh{\widehat V}
\def\Vt{\widetilde V}
\def\Uh{\widehat U}
\def\Ut{\widetilde U}
\def\St{\widetilde S}
\def\Sh{\widehat S}

\def\At{\widetilde A}

\def\Omt{\widetilde\Omega}

\def\Nystrom{Nystr{\"o}m }
\def\Nystromm{Nystr{\"o}m}

\newtheorem{theorem}{Theorem}[section]
\newtheorem{lemma}{Lemma}[section]

\newcommand{\rr}[1]{\textcolor{red}{#1}}

\newcommand{\ignore}[1]{}

\title{
Fast and stable randomized low-rank matrix approximation
}
\author{Yuji Nakatsukasa}

\begin{document}
\maketitle

\begin{abstract}
Randomized SVD has become an extremely successful approach for efficiently computing a low-rank approximation of matrices. 
In particular the paper by Halko, Martinsson, 
and Tropp (SIREV 2011) contains extensive analysis, and has made it a very popular method. 
The typical complexity for a rank-$r$ approximation of $m\times n$ matrices is 
$O(mn\log n+(m+n)r^2)$ for dense matrices. 
The classical \Nystrom method is much faster, but applicable only to positive semidefinite matrices. 
This work studies a generalization of \Nystromm's method applicable to general matrices, and shows that (i) it has near-optimal approximation quality comparable to competing methods, 
(ii) the computational cost is 
the near-optimal $O(mn\log n+r^3)$ for dense matrices, with small hidden constants, and (iii) crucially, it can be implemented in a numerically stable fashion despite the presence of an ill-conditioned pseudoinverse.
 Numerical experiments illustrate that generalized \Nystrom can significantly outperform state-of-the-art methods, especially when $r\gg 1$, achieving up to a 10-fold speedup. 
The method is also well suited to updating and downdating the matrix. 
\end{abstract}





\section{Introduction}
Randomized numerical linear algebra, in particular the randomized SVD by Halko, Martinsson and Tropp~\cite{halko2011finding} has become a highly successful and important practical algorithm for efficiently finding a near-optimal low-rank approximation to a matrix $A\in\mathbb{R}^{m\times n}, m\geq n$. 
In its basic form, the algorithm, which we refer to as HMT (also known as the range finder), 
finds a rank-$r$ ($r\leq \min(m,n)$, usually $r\ll \min(m,n)$) approximant $\hat A_r\approx A$ as follows:
\begin{algorithm}[H]
  \caption{HMT: given 
$A\in\mathbb{R}^{m\times n}$ and $r\in\mathbb{N}$, find a rank-$r$ approximation
$A\approx \hat A_r$. 
}
  \label{alg:HMT}
  \begin{algorithmic}[1]
    \STATE 
Generate a random sketch matrix $\Omega\in\mathbb{R}^{n\times r}$. 
\STATE Compute $A\Omega$. 
\STATE Orthogonalize $A\Omega$ to obtain $Q=\mbox{orth}(A\Omega)$ (e.g. the thin QR factorization $A\Omega=QR$).
\STATE $\hat A_r=Q(Q^T\!A)=(QU_0)\Sigma_0 V_0^T$, where $Q^T\!A=U_0\Sigma_0 V_0^T$ is the SVD. 
  \end{algorithmic}
\end{algorithm}
%

HMT comes with attractive theoretical guarantees~\cite[\S 10]{halko2011finding}, essentially showing that it gives a near-optimal low-rank approximation to $A$ (we make this more precise in Section~\ref{sec:analysis}). 



For the choice of the random sample (or sketch) matrix $\Omega$, structured matrices allowing for fast application such as the
SRHT~\cite{drineas2011faster,boutsidis2013improved,sarlos2006improved}, SRFT~\cite{rokhlin2008fast,tropp2011improved}, and subsampled DCT~\cite{avron2010blendenpik} 
have been proposed. This reduces the sampling cost of forming $A\Omega$ to $O(mn\log n)$. 
Gaussian matrices are however the best understood class of random matrices, with sharp error bounds available for randomized SVD~\cite[\S 10]{halko2011finding}. The overall cost of HMT is $O(mnr)$; this can be reduced to $O(mn\log n +(m+n)r^2)$ by using an interpolative decomposition, at the expense of slightly worse accuracy. 
We recommend the recent survey by Martinsson and Tropp~\cite{MartinssonTroppacta} for an excellent overview of randomized algorithms in numerical linear algebra. 

This paper is about a generalization of \Nystromm's method, 
which is a classical method applicable to positive semidefinite (PSD) matrices  $A\succeq 0$, and finds a rank-$r$ approximation 
\begin{equation}  \label{eq:nystrom}
A\approx AX(X^T\!AX)^\dagger (AX)^T,
\end{equation}
for a sketch matrix $X\in\mathbb{R}^{n\times r}$.
In its original form~\cite{nystrom1930praktische,williams2001using}, \Nystromm's method takes $X$ 
to be a subset of the columns of the identity, so that $AX$ is the corresponding  columns of $A$ and $X^T\!AX$ is $A$'s principal submatrix. 
It has since been generalized to other sketch matrices, and extensively analyzed in \cite{gittens2011spectral,gittens2016revisiting}. 
Unless otherwise mentioned, in this paper we refer to \eqref{eq:nystrom} with a general sketch $X$ as the \Nystrom method. 
 It is a popular method in machine learning for working efficiently with kernel matrices~\cite{fowlkes2004spectral,scholkopf2002learning}. 

The observation that motivated this work is one given in~\cite{halko2011finding,martinsson2016randomized} and explained in \cite{gittens2016revisiting}:
if one takes $X=Q$ from the $Q$ of step 3 of HMT to obtain the approximation 
\begin{equation}  \label{eq:nystHMT}
AQ(Q^T\!AQ)^\dagger (AQ)^T,\quad Q=\mbox{orth}(AX),
\end{equation}
then the resulting accuracy $\|A- AQ(Q^T\!AQ)^\dagger (AQ)^T\|$ is considerably better than $\|A-QQ^TA\|$ with HMT. Note that~\eqref{eq:nystHMT} 
is no more expensive than HMT; it is slightly cheaper. 
Indeed, simple experiments and theory~\cite{gittens2016revisiting} 
reveal that the accuracy of the \Nystrom approximant~\eqref{eq:nystrom} is comparable to (only slightly worse than) HMT 
when one simply takes $X$ to be a random matrix $X=\Omega$ in~\eqref{eq:nystHMT}, which has a much lower cost. 
The complexity of \Nystrom is then the near-optimal $O(n^2\log n +r^3)$. 
To summarize, it appears that for positive semidefinite matrices, 
\begin{itemize}
\item for roughly the same cost, \Nystrom outperforms HMT in accuracy by taking $X=Q$ in~\eqref{eq:nystHMT}, and 
\item for roughly the same accuracy, \Nystrom outperforms HMT in speed significantly by taking the same sketch matrix $X=\Omega$. 
\end{itemize}
For $A\succeq 0$, therefore, \Nystromm's method appears to be the method of choice. 
A noteworthy aspect of \Nystrom with $X$ random (not involving $A$)
 is that the algorithm becomes single-pass, requiring only the linear sketches $AX$ and $X^T\!AX$, thus ideal 
in the streaming model~\cite{tropp2019streaming}. 
The combined efficiency, accuracy and performance of \Nystrom make it a popular algorithm in scientific computing and machine learning. 
However, \Nystrom clearly requires $A\succeq 0$. 
It also involves the (pseudo)inverse of the matrix $(X^T\!AX)^\dagger$, which suggests that numerical instability can be an issue in floating-point arithmetic.

On the other hand, single-pass algorithms have been developed for general nonsymmetric and rectangular $m\times n$ matrices~\cite{tropp2017practical,tropp2019streaming,upadhyay2016fast}. However, unlike \Nystromm, these algorithms require orthogonalization steps, resulting in the complexity $O(mn\log n +(m+n)r^2)$; 
the $(m+n)r^2$ term comes from orthogonalizing an $m\times r$ matrix, which becomes the dominant term when $r\geq \sqrt{\min(m,n)}$. 
Orthogonalization can dominate even when $r$ is much smaller, e.g. in the streaming model~\cite{tropp2019streaming} or when $A$ has additional data-sparse structure. 

In short, the current state-of-the-art 
appears to be that 
\Nystrom is excellent in speed and accuracy for PSD matrices, but for general matrices, one cannot do nearly as well. 
One might wonder if PSD really is a special matrix structure that \Nystrom takes advantage of. 
The above discussion motivates the following questions:

\begin{itemize}
\item Can we generalize the attractive features of \Nystrom (near-optimal accuracy+complexity, storage efficiency, single-pass and no orthogonalization) 
 to general matrices?
\item Is the \Nystromm -like method numerically stable, despite the pseudoinverse? 
\end{itemize}
We answer these questions in the affirmative, and identify an algorithm
which we call \emph{generalized \Nystrom}
that is numerical stable, single-pass, avoids the $O(mr^2)$ orthogonalization cost, and is near-optimal in complexity and accuracy. 
On a desktop machine generalized \Nystrom is seen to outperform HMT by up to an order of magnitude. It is also suitable for updating and downdating the matrix, for which the speedup can be even greater.


\emph{Notation}.
We use matrix norms $\|\cdot\|$  without subscripts for inequalities and arguments that hold for any unitarily invariant norm. We use $\|\cdot\|_F$ for the Frobenius norm and $\|\cdot\|_2$ for the spectral norm. 
$\sigma_i(A)$ denotes the $i$th largest singular value of a matrix $A$. 
For 
matrices $X,Y$ of the same height, 
we use $\mathcal{P}_{X,Y}=X(Y^TX)^\dagger Y^T$ to denote an (oblique) projection onto the column space of $X$. Unless otherwise mentioned we assume that 
$Y^TX$ has full column rank, 
 so the 
 row space of $\mathcal{P}_{X,Y}$ is contained in but not equal to that of $Y^T$. Note that $\mathcal{P}_{XM,YN}=\mathcal{P}_{X,Y}$ for any nonsingular matrices $M,N$. 
$\mathcal{P}_{X}:=\mathcal{P}_{X,X}=X(X^TX)^\dagger X^T$ denotes an orthogonal projection, for which $\|\mathcal{P}_{X}\|_2=1$. 
$\hat A_r$ denotes a rank-$r$ approximant to $A$ with a specified algorithm, and $A_r$ is the rank-$r$ truncated SVD of $A$, which is the optimal rank-$r$ approximant 
in any unitarily invariant norm~\cite[\S 7.4.9]{hornjohn}.
The expected value of a quantity $f(X)$ is denoted by 
$\mathbb{E}f$, where we use subscripts to indicate the random variable as in
$\mathbb{E}_Xf$ when necessary. 
For simplicity we focus on real matrices $A\in\mathbb{R}^{m\times n}$, but everything carries over to complex matrices $A\in\mathbb{C}^{m\times n}$ by replacing the superscript $T$ with $*$, so e.g.\! the HMT approximant becomes $\hat A_r = QQ^*A$. 
\section{The generalized \Nystrom method}\label{sec:gnintro}
How can we ``generalize'' \Nystromm's method~\eqref{eq:nystrom} to nonsymmetric, rectangular matrices? 
Clearly one needs to sketch from the left and right using different matrices; otherwise even the size may not match.
Noting that the \Nystrom approximation~\eqref{eq:nystrom} takes $AX$ and $X^T\!A$ as the column and row spaces respectively, for $A\in\mathbb{R}^{m\times n}$ it is natural to look for 
an approximant with column space $AX$ and row space $Y^T\!A$ for random sketch matrices $X\in\mathbb{R}^{n\times r},Y\in\mathbb{R}^{m\times (r+\ell)}$ ($\ell$  is an oversampling parameter, whose role we discuss later). This leads to an approximant of the form $AXWY^T\!A$, where $W\in\mathbb{R}^{r\times (r+\ell)}$ is a small ``core'' matrix. There are two natural choices of $W$. One is $W=(AX)^\dagger A (Y^T\!A)^\dagger$, which minimizes the Frobenius norm of the error, see e.g.~\cite{cortinovis2019lowrank}, \cite[\S 13]{MartinssonTroppacta}; its computation requires $O(mr^2)$ cost. The other choice $W=(Y^T\!AX)^{\dagger}$ has an interpolatory property 
  and is clearly (cheaper and) closer to the \Nystrom method. Our starting point is therefore the rank-$r$ approximation 
\begin{equation}  \label{eq:start}
A\approx \hat A_r =AX(Y^T\!AX)^{\dagger}Y^T\!A  .
\end{equation}
We refer to this as the \emph{generalized \Nystromm} (GN) method. 
Clearly, it reduces to standard \Nystrom when $A\succeq 0$ and $X=Y$. 
The expression~\eqref{eq:start} is not new; it has been suggested by Clarkson and Woodruff~\cite{clarkson2009numerical,woodruff2014sketching}, with a different derivation based on a Johnson--Lindenstrauss mapping (we discuss their work  more in Section~\ref{sec:related}). 

The approximant~\eqref{eq:start} is very efficient to compute, and yet turns out to have a quasi-optimal approximation guarantee; see \cite{clarkson2009numerical,woodruff2014sketching}, \cite[\S 10]{tropp2017practical} and Section~\ref{sec:gennyst}. 
However, the presence of the matrix (pseudo)inverse  $(Y^T\!AX)^{\dagger}$ is alarming in terms of numerical stability; indeed the matrix $Y^T\!AX$ will almost invariably be ill-conditioned. A naive norm-based stability analysis would bound the error in computing GN by $O(u\kappa_2(Y^T\!AX)^2)$, 
which would mean no accuracy at all in many cases. 
One can understand the work \cite{tropp2017practical} as a stabilized version of~\eqref{eq:start}, by virtue of an orthogonalization step. However, this comes at the cost of an extra $O((m+n)r^2)$ operations. 
As no analysis (and few experiments) appears to have been performed on the numerical stability of~\eqref{eq:start} accounting for roundoff errors in finite-precision arithmetic, and it is unclear whether orthogonalization is really necessary. 

In this work we perform such analysis, and show that, while stability cannot be established for~\eqref{eq:start} as is, there is an inexpensive modification
that guarantees stability: 
\begin{equation}  \label{eq:ournystrom}
A\approx \hat A_r = AX(Y^T\tilde AX)_\epsilon^{\dagger}Y^T\!A,
\end{equation}
which we call the  \emph{stabilized} generalized \Nystrom (SGN) method. 
Here $\tilde A$ is any matrix such that 
$ \tilde A = A+\delta A$ where $\|\delta A\|=O(u\|A\|),$ in which $u$ is the unit roundoff ($u\approx 10^{-16}$ in standard IEEE double precision), and $(Y^T\!AX)_\epsilon^{\dagger}$ denotes the $\epsilon-$pseudoinverse, that is, 
if $Y^T\!AX=[U_1,U_2]
\big[
\begin{smallmatrix}
\Sigma_1 & \\ & \Sigma_2  
\end{smallmatrix}
\big]
 [V_1,V_2]^T$ is the SVD where $\Sigma_1$ contains singular values larger than $\epsilon$, then $(Y^T\!AX)_\epsilon^{\dagger}=V_1\Sigma_1^{-1}U_1^T$. 
In this paper we always take $\epsilon=O(u\|A\|)$, a modest multiple of the unit roundoff $u$ times $\|A\|$. 

A careful inspection reveals that the expression~\eqref{eq:ournystrom} is somewhat redundant: any $AX(Y^T\tilde AX)_\epsilon^{\dagger}Y^T\!A$ is equal to 
$AX(Y^T\tilde A_2X)^{\dagger}Y^T\!A$ where $\tilde A_2$ is also within $\epsilon$ of $A$. 
In an actual computation, one would always attempt to evaluate $AX(Y^TAX)_\epsilon^{\dagger}Y^T\!A$; 
the point of SGN~\eqref{eq:ournystrom} is that a (carefully) computed 
approximation $fl(AX(Y^TAX)_\epsilon^{\dagger}Y^T\!A)$
can be written exactly 
in the form~\eqref{eq:ournystrom} for some $\tilde A$ (row-wise; see Section~\ref{sec:gennyst}), 
 and any approximant of the form~\eqref{eq:ournystrom} has error $\|A-fl(\hat A_r)\|$ comparable to \eqref{eq:start}; SGN is therefore a stable method. A proof of this claim is a key contribution of this paper. 



Many properties (complexity, performance and even stability in practice) are shared between generalized \Nystromm~\eqref{eq:start} and the stabilized version~\eqref{eq:ournystrom}, so in what follows, when we simply refer to GN, the arguments apply both to~\eqref{eq:start} and~\eqref{eq:ournystrom} unless otherwise mentioned. When the distinction is important we call \eqref{eq:ournystrom} \emph{stabilized} GN, and \eqref{eq:start} \emph{plain} GN. 
It turns out that, as we highlight in Section~\ref{sec:gennyst}, while 
stability cannot be established for plain GN \eqref{eq:start}, 
its instability is benign, and one usually obtains satisfactory results. 

Below is a pseudocode for plain and stabilized GN. 

\begin{algorithm}[H]
  \caption{(stabilized) Generalized \Nystrom: given 
$A\in\mathbb{R}^{m\times n}$ and $r$, find a rank-$r$ approximation
$A\approx \hat A_r$. 
}
  \label{alg:GN}
  \begin{algorithmic}[1]
    \STATE Generate sketch matrices $X\in\mathbb{R}^{n\times r}$, $Y\in\mathbb{R}^{m\times (r+\ell)}$, where $0<\ell=\left \lceil{0.5r}\right \rceil $ is suggested. 
    \STATE Compute $AX$, $Y^T\!A$, and $QR$ factorization $Y^T\!AX=QR$. 
\STATE $\hat A_r= ((AX) R^{-1})(Q^T(Y^T\!A))$, or \hfill (plain GN, stable most of the time) \\
 $\hat A_r= ((AX) R_\epsilon^\dagger)(Q^T(Y^T\!A))$ 
\hfill (stabilized GN)
  \end{algorithmic}
\end{algorithm}
The outputs of (S)GN are $AX\in\mathbb{R}^{m\times r}$, $Y^T\!A\in\mathbb{R}^{(r+\ell)\times n}$, and the small matrices $R\in\mathbb{R}^{r\times r}$ and $Q\in\mathbb{R}^{(r+\ell)\times r}$. The memory requirement is $\approx mr+n(r+\ell)+r^2\approx (m+1.5n)r$ with the recommended choice $\ell=0.5r$. 


While $\hat A_r = ((AX)R^{-1})(Q^T(Y^T\!A))$, 
evaluating $(AX)R^{-1}$ or $Q^T(Y^T\!A)$ is usually not advisable, as that would require $O(mr^2)$ cost. 
Instead, one would use the output factors to perform further operations with $\hat A_r$. 
For example, to compute a matrix product
$AW$ for a given $W\in\mathbb{R}^{n\times k}$, one can perform 
$AW \approx \hat A_rW= AX(R^{-1}(Q^T(Y^T\!A W)))$
in the order indicated. This requires $O((m+n)rk)$ operations, the same complexity required by performing $WA$ with other methods such as HMT.
Similarly, for left-multiplication $WA$ one would perform
$WA\approx W\hat A_r = (((WAX)R^{-1})Q^T)Y^T\!A$. 



Here the ($\epsilon$-)pseudoinverses $(Y^T\!AX)^{\dagger}, (Y^T\!AX)_\epsilon^{\dagger}$ are implemented via a QR factorization; one could also use the SVD. 
We discuss implementation details, including computing $R_\epsilon^\dagger$ and 
the choice of $\ell$, in Section~\ref{sec:implement}. 

\subsection{Numerical illustration}
To motivate the study of GN, let us illustrate its performance in comparison with popular methods when applied to positive definite matrices. 
We generate a $50000\times 50000$ positive definite matrix $A=Q\Lambda Q^T$, where $Q$ is a random orthogonal matrix (Q-factor in the QR factorization of a square Gaussian matrix) and $\Lambda$ has geometrically decaying singular(=eigen) values. 
We compare the speed and accuracy $\|A-\hat A_r\|_F$, where $\hat A_r$ is a rank-$r$ approximant, varying $r$ from $10^3$ to $10^4$. $X,Y$ are subsampled DCT matrices. 

\begin{figure}[htpb]
  \begin{minipage}[t]{0.5\hsize}
      \includegraphics[height=55mm]{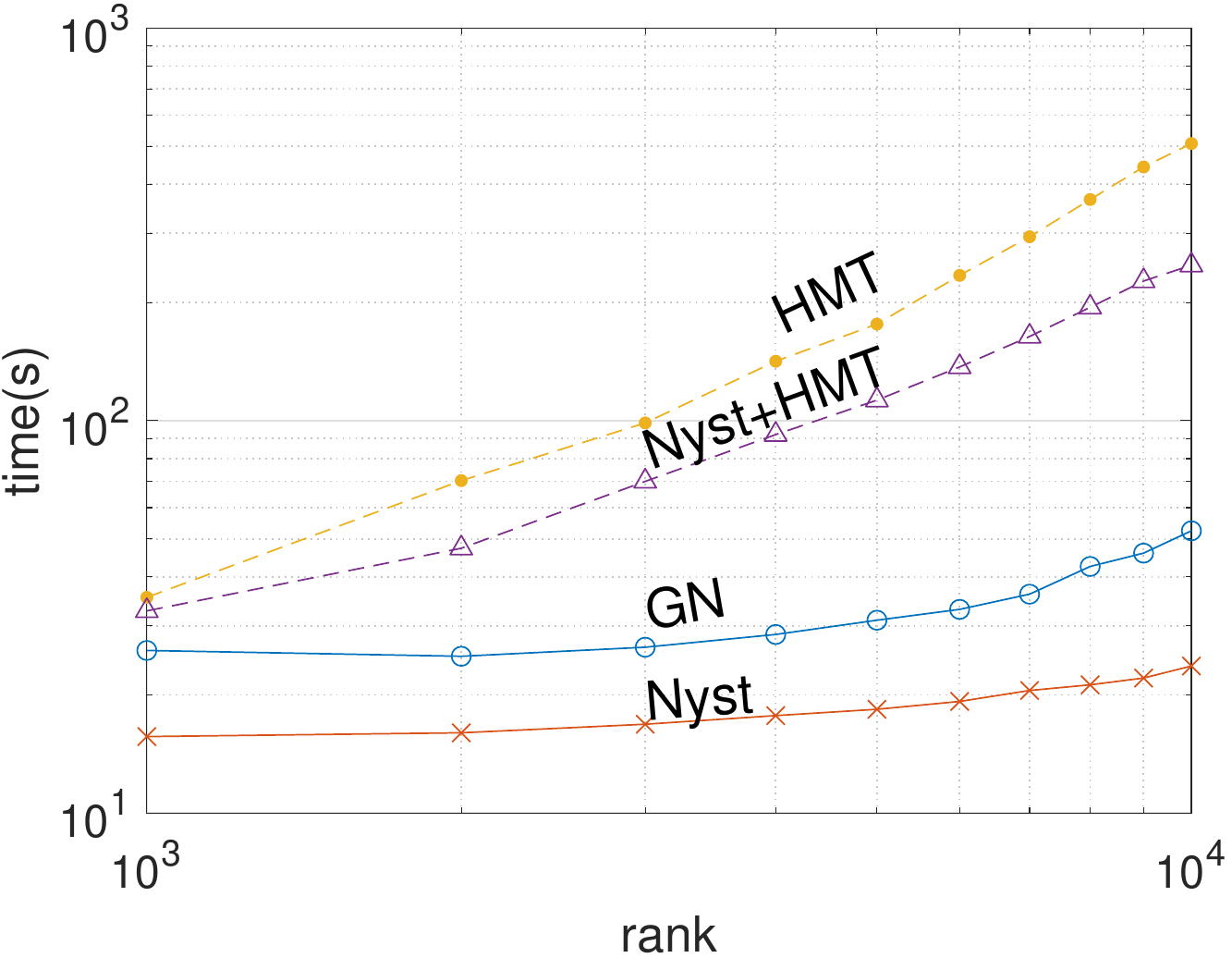}      
  \end{minipage}   
  \begin{minipage}[t]{0.5\hsize}
      \includegraphics[height=55mm]{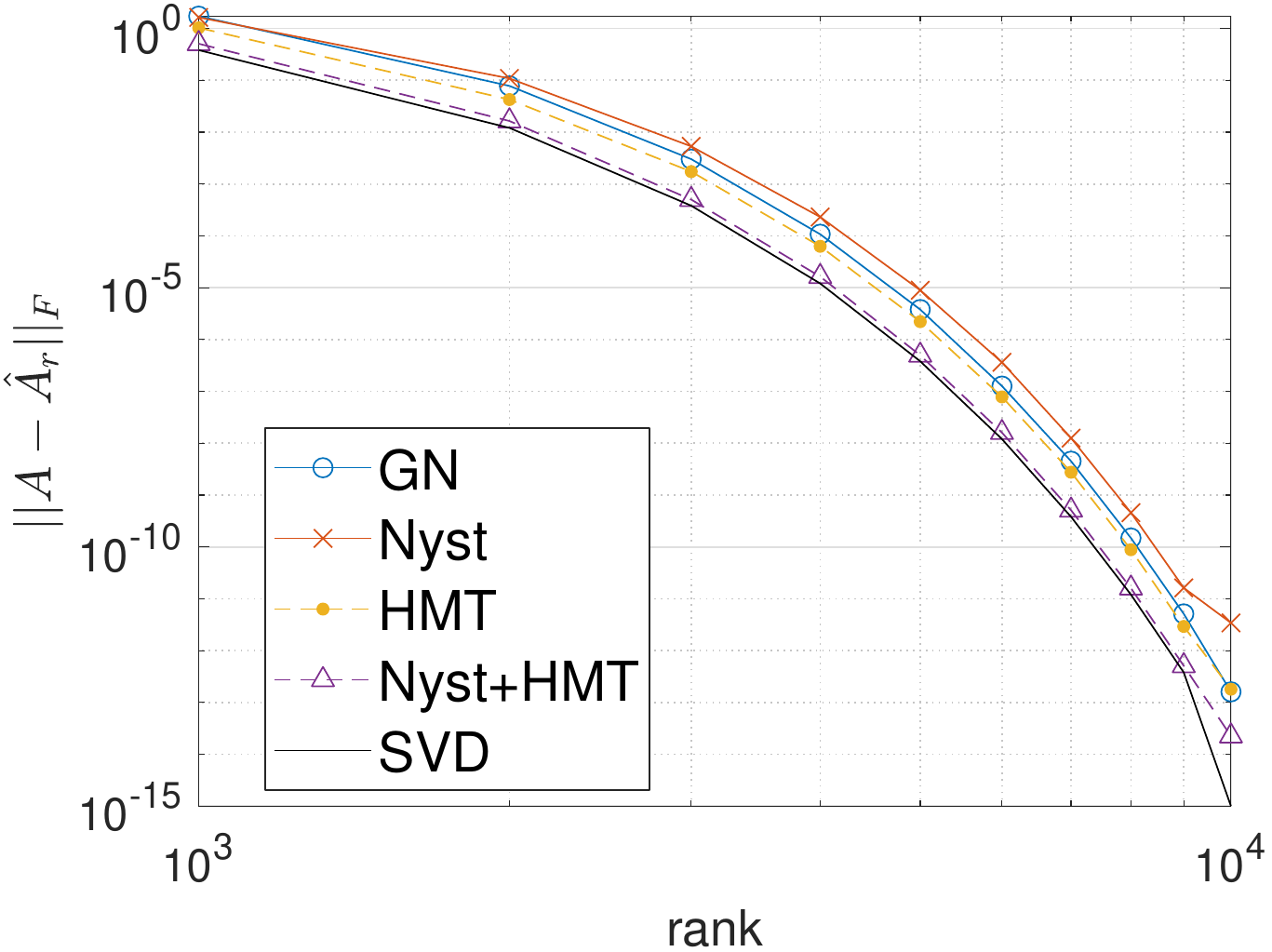}            
  \end{minipage}
  \caption{Algorithms for $A\succeq 0$.
GN is the generalized \Nystrom method~\eqref{eq:start}, 
Nyst is \Nystromm's method, HMT is Halko-Martinsson-Tropp~\cite{halko2011finding} and HMT+Nyst is~\eqref{eq:nystHMT}. SVD is the optimal truncated SVD. 
}
  \label{fig:posdef}
\end{figure}

The results are shown in Figure~\ref{fig:posdef}. 
The \Nystrom method~\eqref{eq:nystrom} (shown as Nyst) runs the fastest, about 20x faster\footnote{Throughout, numerical experiments were performed in MATLAB version 2020a on a desktop computer with 256GB memory. The runtime of HMT includes the small SVD $Q^T\!A=U_0\Sigma_0 V_0^T$; without it, HMT has about the same speed as Nyst+HMT.}  than HMT (Algorithm~\ref{alg:HMT})
when $r=10^4$.  
GN \eqref{eq:ournystrom} (the plain~\eqref{eq:start} actually performs very similarly, see Section~\ref{sec:implement})
has runtime roughly twice that of \Nystromm, 
 and up to 10x faster than HMT, with larger speedup observed for larger rank $r$, reflecting the lower complexity when the $O(mr^2)$ orthogonalization cost is dominant. 

In terms of accuracy, the main message is that all methods are close to optimal, tracking the optimal truncated SVD (shown as SVD) 
to within a modest factor. 
Nyst+HMT implements~\eqref{eq:nystHMT}, and comes the closest to optimal. 
Importantly, no numerical instability is observed in any method, despite the presence of the pseudoinverses in \Nystrom and GN. The accuracy of GN is only marginally worse than HMT, a difference that is unlikely to matter when the singular values decay sufficiently fast. 

As is well known with randomized SVD methods~\cite[\S 7]{halko2011finding}, all methods are remarkably consistent: Despite the random nature, running the experiment multiple times results in nearly identical figures. 

All methods compared here, except classical \Nystromm, are applicable to general matrices. Overall, generalized \Nystrom is seen to extend the attractive properties (speed+accuracy) of \Nystrom to nonsymmetric and rectangular matrices. 
In the remainder of this paper we study the generalized \Nystrom method in detail. 


\subsection{Properties of generalized \Nystromm}
Generalized \Nystrom has the following properties: 
\begin{itemize}
\item It is based on linear sketches: $AX$, $Y^T\!A$, and $Y^T\!AX$ (and no terms of the form e.g. $AA^T,A^T\!A$,\ldots). 
\item Its cost is $O(N_r+r^3)$, where $N_r$ is the cost for forming $AX$ and $Y^TA$, where $\ell=O(r)$ is assumed. 
Specifically, the cost is $O(mn\log r+r^3)$ if $A$ is dense (using the SRFT sketch matrices), and $O(\mbox{nnz}(A)r+r^3)$ if $A$ is a sparse matrix with 
$ \mbox{nnz}(A)$ nonzero elements. 
\item The approximation quality $\|A-\hat A_r\|$ is near optimal, on the order of the error $\|A- A_{\hat r}\|$ with the truncated SVD for some $\hat r$ slightly smaller than $r$.
\item 
The stabilized version~\eqref{eq:ournystrom}
can be implemented in a numerically stable manner in the presence of roundoff errors. 
\end{itemize}
The first two points are straightforward to verify, and they make GN among the most efficient methods for computing low-rank approximations for general matrices. 

The third and fourth points are the main technical results of this paper, and 
treated in Sections~\ref{sec:analysis} and \ref{sec:gennyst}. 
To guarantee numerical stability a careful implementation is required, as we discuss in Section~\ref{sec:implement}. 
%



It is worth noting that the output of GN does not give an (approximate) truncated SVD, as no factor has orthonormal columns. 
This is expected of an algorithm that requires less than $O(\max(m,n)r^2)$ operations, 
which would be needed for simply orthogonalizing a matrix of size $\max(m,n)\times r$. 
This does come with limitations: we do not have access to information available with approximants with orthogonal factors, such as the singular values and vectors of the approximant $\hat A_r$, which can be used as approximate singular values of those of $A$. 
For example in HMT, $(QU_0)\Sigma_0 V_0^T\approx A$ is an approximate SVD. 


\subsection{Near-optimal complexity}\label{sec:cost}
Generalized \Nystrom is essentially optimal in computational complexity 
in a number of natural settings. 
First when $A$ is dense, clearly at least $O(mn)$ operations are necessary for a low-rank approximation, as failure to read one element can result in poor approximation. 
With a fast sketching matrix such as SRFT, 
computing $AX,Y^T\!A$  can be done in $O(mn\log n)$ operations.  
Thus generalized \Nystrom is $O(mn\log n+r^3)$ flops, which can be further improved to $O(mn\log r+r^3)$~\cite{woolfe2008fast}. The first term is clearly optimal up to $O(\log n)$ or $O(\log r)$, and $O(r^3)$ operations is the standard cost for working with $r\times r$ full-rank unstructured matrices (with a Strassen-like fast matrix multiplication algorithm one could reduce it to $r^{\omega}$ where currently $\omega\approx 2.4$). 

For generalized \Nystromm, it is even possible to work out the hidden constants, which are not large. 
For the SRFT, the sketching requires $N_r=10mn\log n$ flops, since the FFT with an $n\times n$ matrix costs $5n\log n$ flops per vector~\cite[\S 4.6]{Golubbook2}. 
For the core matrix $Y^T\!AX$, computing $(Y^T\!A)X$ is just $5nr\log n$  flops, and the pseudoinverse based on the QR factorization (see Section~\ref{sec:pinv}) requires $2(r+\ell)r^2-\frac{2}{3}r^3$~\cite[\S 5.2.9]{Golubbook2}, which is $\frac{7}{3}r^3$ when $\ell=r/2$ (our recommended choice).
The overall cost of generalized \Nystrom is thus $10mn\log n+\frac{7}{3}r^3$ flops (plus strictly lower order terms, such as $5nr\log n$ and $O(r^2)$). 

Similar near-optimality can be established in the streaming model~\cite{tropp2017practical,tropp2019streaming}, in which updates of $A$ are given sequentially and past information is discarded. 

The situation can be different when $A$ is structured, so that multiplying $A$ by a vector can be done efficiently. An example is when $A$ sparse, as $mn\log n$ may be much larger than $\mbox{nnz}(A)$. In this case, one can take the sketch matrices to be e.g.\! the CountSketch matrix with one nonzero element of $\pm 1$ in each row, achieving input-sparsity time~\cite{clarkson2017low}. Such choice reduces the cost in computing $AX, Y^T\!A$ to $O(\mbox{nnz}(A))$. 
However, the analysis in \cite{clarkson2017low} suggests that one would need to take $\ell=O(r^2)$ to guarantee good accuracy.

In the sparse (and more generally structured) case one might naturally require that the factors preserve the sparsity, a property that GN may not satisfy. Alternatives such as the CUR factorization~\cite{goreinov1997theory,mahoney2009cur} may be more attractive in such settings. 


\subsection{Related studies}\label{sec:related}
While we have derived the GN approximant \eqref{eq:ournystrom} by generalizing the \Nystrom method, 
the expression~\eqref{eq:ournystrom} itself is not new; it appears in Clarkson and Woodruff~\cite[Thm.~4.7]{clarkson2009numerical} and \cite[Thm 4.3]{woodruff2014sketching}. However, as mentioned above, 
the treatment there focuses on the case where $X,Y$ are Countsketch or sign matrices, and 
their results suggest that a significant oversampling $\ell\gg r$ would be necessary for near-optimal accuracy. 
Also, the fact that the $O(mr^2)$ cost due to orthogonalization can be avoided is not emphasized in their work, with no experiment reported to illustrate the practical speed. 
Finally, and most significantly, the crucial issue of numerical stability in finite-precision arithmetic is not discussed.

Indeed, avoiding the potential instability is the primary reason the authors in \cite{tropp2017practical} orthogonalize $AX=QR$ and take $A\approx Q(Y^TQ)^{\dagger}Y^T\!A$, which is mathematically identical to GN~\eqref{eq:start}. 
However, while \cite{tropp2017practical} thus avoids the numerical instability
 resulting from inversions\footnote{To be precise, however, we cannot immediately conclude that \cite{tropp2017practical} is numerically stable: The QR factorization $AX=QR$ in finite-precision arithmetic results in a computed $Q$ that has error $O(u\kappa_2(AX))$~\cite[Ch.~19]{Higham:2002:ASNA}, so strictly speaking the stability of~\cite{tropp2017practical}---or even HMT---is an open problem. Of course numerical evidence strongly suggest they are indeed stable. 
In any case our focus is generalized \Nystromm~\eqref{eq:ournystrom} and the proof of its numerical stability when implemented appropriately. 
}, orthogonalization requires at least $O(mr^2)$ operations, which is often the dominant part of the cost~\cite{tropp2019streaming}.  
Similarly, 
HMT is also based on orthogonal projection (and hence ``stable''), but inevitably involves the $O(mr^2)$ orthogonalization cost.

The equivalence between plain GN~\eqref{eq:start}
and \cite{tropp2017practical} means that the extensive accuracy analysis there (in exact arithmetic) carries over verbatim, some of which we rederive in Section~\ref{sec:analysis}. 
Furthermore, we carefully investigate the effects of roundoff errors in finite-precision arithmetic, to show that with an appropriate implementation, the presence of the pseudoinverse in GN is not a problem for stability, even though the core matrix $Y^T\!AX$ does become severely ill-conditioned. 




\subsubsection{Comparison with other algorithms}\label{sec:compother}

It turns out that many algorithms can be written formally as~\eqref{eq:start} for certain choices of $X$ and $Y$.
A key step is to note that 
$ (AX(Y^T\!AX)^{\dagger}Y^T)A  = \mathcal{P}_{AX,Y}A$, where 
$\mathcal{P}_{AX,Y}:=AX(Y^T\!AX)^{\dagger}Y^T$ is an oblique projection onto the column space of $AX$, 
and also $A(X(Y^T\!AX)^{\dagger}Y^T\!A) = A\mathcal{P}_{X,A^TY}$, 
where $\mathcal{P}_{X,A^TY}:=X(Y^T\!AX)^{\dagger}Y^T\!A$ is an oblique projection onto a subspace of the rows of $Y^T\!A$. 
That is, \eqref{eq:start} can be thought of as a two-sided projection of $A$. Indeed we have 
$\mathcal{P}_{AX,Y}A=A\mathcal{P}_{A^TY}=\mathcal{P}_{AX,Y}A\mathcal{P}_{A^TY}$, as can be verified by 
\[
\mathcal{P}_{AX,Y}A\mathcal{P}_{X,A^TY} = AX(Y^T\!AX)^{\dagger}Y^T\!AX(Y^T\!AX)^{\dagger}Y^T\!A=AX(Y^T\!AX)^{\dagger}Y^T\!A. 
\]
We now show that many methods can be mathematically (ignoring roundoff errors) written as~\eqref{eq:start}. 
First, clearly \Nystrom is \eqref{eq:start} with $X=Y$. 
Clarkson and Woodruff~\cite{clarkson2009numerical} is the plain generalized \Nystrom where $X,Y$ are Countsketch matrices. Its equivalence (up to the choice of $X,Y$) to Tropp et al.~\cite{tropp2017practical} is explained in \cite{tropp2017practical} and mentioned above. 

We now show that the HMT approximant 
$QQ^T\!A$ (where $AX=QR$ is the QR factorization) is~\eqref{eq:start} with $Y=AX$. To see this, note that 
$AX((AX)^T\!AX)^{\dagger}(AX)^T\!A=P_{AX}A=QQ^T\!A$. 
This choice results in the advantage that the projection $\mathcal{P}_{AX}$ is orthogonal, so $\|\mathcal{P}_{AX}\|_2=1$. 

One can obtain improved approximants by taking higher powers, for example 
 $A(A^T\!A)^pX(Y^T\!A(A^T\!A)^pX)^{\dagger}Y^T\!A$ where $p\geq 0$, which corresponds to subspace iteration extensively studied by Gu~\cite{gu2015subspace}. 
In this case we have an orthogonal projection by taking $Y=(A^T\!A)^pX$, since then the approximant can be written 
$A(A^T\!A)^pX(X^T(A^T\!A)^pA(A^T\!A)^pX)^{\dagger}X^T(A^T\!A)^pA=Y(Y^TY)^\dagger Y^T\!A$. For numerical stability, it is advisable to compute the thin QR factorization of $Y=QR$ and compute $QQ^T\!A$. 

We summarize and compare these algorithms in Table~\ref{tab:tab}. 

\begin{table}[htbp]
  \centering
  \caption{Comparison of randomized algorithms for rank-$r$ approximation $A\approx AX(Y^T\!AX)^{\dagger}Y^TA \in\mathbb{R}^{m\times n}$. 
$\Omega,\tilde\Omega\in\mathbb{R}^{m\times O(r)}$ represent random sketch matrices. 
$N_r$ (and $\tilde N_r$, see text below) is the cost of forming the products 
$A\Omega,\tilde\Omega^TA$.  
($\surd$) means the method is observed and conjectured to be numerically stable in practice but no proof is available. ($\times$) means the method is unstable, but usually performs in a stable way, see Section~\ref{sec:gennyst}. 
Here we assume $m\geq n$ to simplify the complexity. 
}
  \label{tab:tab}
\small
  \begin{tabular}{c|c|c|c|c}
    & $X,Y$ & $p$  & Stable? & cost for dense $A$\\\hline
HMT~\cite{halko2011finding} & $X=\Omega, Y=AX$ & 0 & ($\surd$) & $O(\tilde  N_r+mr^2)$  \\
\Nystrom($A\succ 0$)~\cite{gittens2016revisiting} & $Y=X=\Omega$ & 0  & $\times$& $O(N_r+mr^2)$ \\
HMT+\Nystrom($A\succ 0$)~\cite{halko2011finding} & $Y=X=Q, A\Omega=QR$& 1 & $\times $& $O(N_r+mr^2)$ \\
Subspace iteration~\cite{gu2015subspace}& $X=(A^T\!A)^p\Omega, Y = AX$& $p>1$ & $(\surd)$& $O(2p\tilde N_r+mr^2)$ \\
Demmel-Grigori-Rusciano~\cite{demmelgrigori2019improved}& \cite{clarkson2009numerical}+extra term & 0 & ($\times$) & $O(N_r+mr^2)$ \\
Tropp17~\cite{tropp2017practical} & $X=\Omega,Y=\tilde\Omega$ & 0 & ($\surd$) & $O(N_r+mr^2)$ \\
plain GN=Clarkson-Woodruff~\cite{clarkson2009numerical} & $X=\Omega,Y=\tilde\Omega$ & 0 & ($\times$) & $O(N_r+r^3)$ \\
stabilized GN~\eqref{eq:ournystrom} & $X=\Omega,Y=\tilde\Omega$ & 0 & $\surd$ & $O(N_r+r^3)$ \\
  \end{tabular}
\normalsize 
\end{table}

The term $O(N_r)$ in the cost is for computing the sketches $AX,Y^T\!A$ ($r$ is essentially the width of $X,Y$), and the specific value depends on the situation. For example for dense matrices, one can use the SRFT matrices for $X,Y$ allowing for fast matrix-vector multiplications $N_r=10mn\log n$ (note its independence of $r$), resulting in the overall cost $O(mn\log n+r^3)$. This applies to all algorithmns but subspace iteration and HMT, which require revisiting the matrix when computing $Q^T\!A$, wherein fast sampling is no longer available. This is why the cost in the table is shown as $\tilde N_r$, which becomes $\tilde N_r=O(mnr)$ in the dense case. 

The bottom three algorithms all start from the same mathematical expression~\eqref{eq:nystrom}, but have important differences. Namely, Tropp17 and Clarkson--Woodruff (which is essentially generalized \Nystromm~\eqref{eq:start}) are mathematically the same, but Tropp17 uses orthogonalization to obtain a (empirically---indicated by $(\surd)$ in the table) stable algorithm, at the cost of the additional $O(mr^2)$ operations. SGN stabilizes Clarkson--Woodruff by a simpler remedy without requiring orthogonalization. 

The recent preprint by Demmel, Grigori and Rusciano~\cite{demmelgrigori2019improved} studies low-rank approximants from the perspective of the LU factorization, and derives an approximant that has an extra term 
of the form $(Y^T)^\dagger M$, where $M\in\mathbb{R}^{(r+\ell)\times n}$ in addition to \eqref{eq:start}. 
The authors show its accuracy is between that of GN and HMT; experiments suggest the accuracy improvement over GN is usually marginal, which is perhaps expected as the random matrix $(Y^T)^\dagger$ may have nothing to do with the column space of $A$. The method is also more expensive than generalized \Nystromm. 

It is worth noting that when $r=O(1)$, the orthogonalization cost $O(mr^2)$ is dominated by the sketch $N_r=O(mn\log n)$ (or $N_r=O(mnr)=O(mn)$, by taking $X,Y$ to be Gaussian), so all the algorithms in Table~\ref{tab:tab} would be optimal; this applies also to a (block-)Krylov subspace method. The advantage of generalized \Nystrom becomes significant when $r\gg 1$, 
as we illustrate in our experiments in Section~\ref{sec:ex}.

An extreme case of GN is when we take $X,Y$ to be subsets of the identity matrix, so that $AX,Y^T\!A$ are simply $A$'s selected columns and rows. For $A\succeq 0$ this is the classical \Nystrom method (with $X=Y$), which has an extremely low $O(r^3)$ complexity. With a random choice of such $X,Y$, this is guaranteed to succeed provided that the matrix is incoherent, as studied in detail by Chiu and Demanet \cite{chiu2013sublinear}. 
They also mention the use of subsampled Fourier matrices to reduce the coherence to obtain an algorithm that works for any matrix; however, their algorithm does not allow for oversampling $\ell>0$, which is crucial for stability as we shall see. 
When $X,Y$ are (carefully chosen) subsets of columns of the identity, the resulting approximant is a so-called CUR factorization  $A\approx CUR$~\cite{goreinov1997theory,mahoney2009cur} 
 in which $C,R$ are subsets of the columns and rows of $A$. The paper \cite{goreinov1997theory} proves existence of a CUR factorization with approximation quality optimal up to a factor $O(\sqrt{k}(\sqrt{m}+\sqrt{n}))$; this takes $U=C^\dagger A R^\dagger$. 
When the core matrix is taken to be $(Y^T\!AX)^{\dagger}$, a recent result \cite{zamarashkin2018existence} proves optimality can be achieved up to a polynomial factor in $r$, and \cite{cortinovis2019lowrank} derives a practical algorithm for finding it. However, the algorithm starts with the SVD of $A$, so it is clearly not competitive with generalized \Nystrom in efficiency.

Other methods based on linear sketching that are not in the form 
$AX(Y^T\!AX)^{\dagger}Y^T\!A$ include 
Boutsidis et al.~\cite{boutsidis2016optimal},
Cohen et al.~\cite{cohen2015dimensionality}
and 
Tropp et al.~\cite{tropp2019streaming}. 
These methods all involve orthogonalization, so require at least $O(N_r+mr^2)$ operations. For positive definite matrices, a sublinear-time algorithm
 for low-rank approximation has been developed recently~\cite{musco2017sublinear}. 


\subsection{Analyzing a numerical algorithm}
To prove that a numerical algorithm is ``good'', 
two facts need to be established: 
\begin{enumerate}
\item[(i)] The algorithm outputs a good approximate solution in exact arithmetic: in our case $\|A-\hat A_r\|$ decays as $r$ increases if $A$ has rapidly decaying singular values. 
\item[(ii)] Roundoff errors do not impair the quality of the output, that is, $\|\hat A_r-fl(\hat A_r)\|$ is small (here $fl(\cdot)$ denotes a quantity computed in a standard IEEE floating-point arithmetic), so that the computed output is still a good approximation. 
\end{enumerate}
Most studies on randomized SVD so far have naturally focused on (i), giving estimates for the quality of the low-rank approximants obtained by randomized SVD algorithms. It appears that little has been done to study the effect of roundoff errors for randomized algorithms. 
Regularization methods have been introduced~\cite{gittens2013topics,li2017algorithm,tropp2017fixed,tropp2017practical}, 
 based on the principle that in order to obtain a numerically stable method, one must avoid inverting ill-conditioned matrices. 
For example, the method in~\cite{tropp2017fixed} works with a shifted matrix $A+\sigma I$ to improve conditioning for $A\succeq 0$; but this technique is conveniently possible only for $A\succeq 0$. 
We shall show that such regularization is often overkill. 



In the following we deal with both (i) approximation accuracy (in Section~\ref{sec:analysis}), and (ii) numerical stability of GN (in Section~\ref{sec:gennyst}).

\section{Approximation accuracy}\label{sec:analysis}
Here we analyze the approximation accuracy (or error) $\|A-\hat A_r\|$, ignoring the effect of roundoff errors. 
While we focus on the HMT and generalized \Nystrom methods, our analysis can be applied to any method based on projection, which includes all algorithms in Table~\ref{tab:tab} but~\cite{demmelgrigori2019improved}. 
 While results on HMT and (plain) GN can be found in the literature, we believe the analysis here is simpler than most, and treats many methods in a unified fashion. 
In addition, based on the analysis we reexamine the oversampling parameter $\ell$ in Section~\ref{sec:implement}, and argue that $r=0.5$ is a safe choice in all cases. This choice is simpler than the recommendations given in \cite{tropp2017practical}.

\subsection{Key facts}
The following facts will be used repeatedly in the forthcoming analysis. 
\begin{itemize}
\item ``Rectangular random matrices are well-conditioned''. 
This informal statement has been made very precise in the literature:
\begin{itemize}
\item For random matrices with i.i.d entries (including Gaussian matrices) of zero mean and unit variance, the classical Marchenko-Pastur (M-P) rule~\cite{pastur1967distribution,yaskov2016short} shows that an $m\times n$ ($m\geq n$) 
matrix $G_{m,n}$ 
has singular values supported 
in the interval $[\sqrt{m}-\sqrt{n},\sqrt{m}+\sqrt{n}]$ (the M-P rule identifies the precise limiting distribution). Extensions have been derived
by Davidson and Szarek~\cite{davidson2001local}, who show for Gaussian matrices that the probability that a singular value lies $\delta$-outside of this interval decays squared-exponentially with $n$ and $\delta$. More precisely, they show that \cite[Thm. II.13]{davidson2001local} (see also \cite[Remark 8.5]{MartinssonTroppacta})
\begin{equation}  \label{eq:Davidson}
\max(\mathbb{P}[\sigma_{n}(G_{m,n})\leq \sqrt{m}-\sqrt{n}-t],
\mathbb{P}[\sigma_{1}(G_{m,n})\geq \sqrt{m}+\sqrt{n}+t])
<\exp(-\frac{n}{2m}t^2). 
\end{equation}
More precisely, the extremal singular values are known to follow the Tracy-Widom distribution~\cite{feldheim2010universality}. 
We also mention Rudelson and Vershynin~\cite{rudelson2009smallest}, who argue that the smallest singular value behaves approximately like $\sqrt{m}-\sqrt{n-1}$. 
This means that the condition number is essentially bounded by $\frac{\sqrt{m}+\sqrt{n}}{\sqrt{m}-\sqrt{n-1}}$, which is modest
if $m\neq n$, and tends to $1$ as $m/n\rightarrow 1$. 
A corollary is that the pseudoinverse typically has norm $\|X^\dagger\|\approx 1/(\sqrt{m}-\sqrt{n-1})$, with high concentration if $m\neq n$. 
\item The result has been qualitatively extended 
to other structured matrices including the SRHT~\cite{drineas2011faster,boutsidis2013improved,sarlos2006improved}, SRFT~\cite{rokhlin2008fast,tropp2011improved}, subsampled DCT matrices~\cite{avron2010blendenpik}, and the CountSketch~\cite{clarkson2017low} matrices. These results are often derived based on the celebrated Johnson-Linderstrauss lemma. The bounds are however usually not as sharp as~\eqref{eq:Davidson}. 
\end{itemize}
\item Properties of projection. 
  \begin{enumerate}
  \item[(i)] Projections annihilate a certain subspace: in particular 
for $\mathcal{P}_{AX,Y}=AX(Y^T\!AX)^\dagger Y^T$ with 
$Y^T\!AX$ having full column rank, 
\begin{equation}  \label{eq:anni}
(I-\mathcal{P}_{AX,Y})AX=0  .
\end{equation}
  \item[(ii)] Norm identity 
    \begin{equation}\label{eq:projs}
\|I-P\|_2=\|P\|_2      , 
    \end{equation}
 which holds for any projection $P$ s.t. $P^2=P$. See Szyld~\cite{szyld2006many} for a delightful account of this important identity. 
  \end{enumerate}
\end{itemize}

Let us note in passing that the fact that rectangular random matrices are well-conditioned is at the heart of a number of recent developments in computational mathematics; 
the randomized least-squares solver Blendenpik~\cite{avron2010blendenpik}, 
the RIP in compressed sensing~\cite{foucart2013mathematical}, 
and stable approximation of functions via discrete least-squares~\cite{cohen2013stability} 
 are among the representative examples. 

Below we sometimes restrict the analysis to Gaussian sketch matrices $X,Y$, when they allow for precise bounds; we will make it clear when this assumption is made. A convenient fact about Gaussians is orthogonal invariance: if $G$ is Gaussian, then so are $GQ_1$ and $Q_2G$ for any orthogonal $Q_1,Q_2$. 
In most cases, the same qualitative results hold for other classes of $X,Y$
with different constants; see e.g.~\cite[Thm.~11.2]{halko2011finding} and \cite[Remark 8.2]{MartinssonTroppacta}. Moreover, while for brevity we focus on bounding the expected value of the error, 
probabilistic error bounds that are satisfied with high probability also hold for GN, analogous to those known for HMT~\cite[\S 10.3]{halko2011finding}. 

\subsection{Key steps for analyzing $\hat A_r-A$}
Recall from Section~\ref{sec:compother} that we can write $\hat A_r=AX(Y^T\!AX)^{\dagger}Y^TA$
in both HMT and (plain) generalized \Nystromm. 
It follows that the error 
 can be expressed in three equivalent forms as
\begin{align}\nonumber 
E:=A-X(Y^T\!AX)^{\dagger}Y^T\!A& =   (I-\mathcal{P}_{AX,Y})A = A(I-\mathcal{P}_{X,A^TY}) = (I-\mathcal{P}_{AX,Y})A (I-\mathcal{P}_{X,A^TY}) . 
\end{align}
In view of the second key fact~\eqref{eq:anni} above, it follows that 
\begin{equation}
  \label{eq:Eform}
  \begin{split}
E &= (I-\mathcal{P}_{AX,Y})A = (I-\mathcal{P}_{AX,Y})A(I-XM_X^T)
  \end{split}
\end{equation}
for any matrix $M_X\in\mathbb{R}^{r \times n}$. 
Below we shall take specific choices of $M_X$ to analyze $\|E\|$ for different algorithms. 

\subsection{HMT}
As discussed in Section~\ref{sec:compother}, HMT (Algorithm~\ref{alg:HMT}), $Y$ is chosen as  $Y=AX$ and thus $\mathcal{P}_{AX,Y}=\mathcal{P}_{AX}$ is an orthogonal projector, so $\|\mathcal{P}_{AX}\|_2=\|I-\mathcal{P}_{AX}\|_2=1$. 

Thus from the first expression in~\eqref{eq:Eform} we obtain\footnote{The analysis here is inspired by an observation in~\cite[\S~4]{sorensen2016deim}, which is attributed to Ipsen.
Other (lengthier) derivations in the literature include the original \cite{halko2011finding}, and those based on the Schur complement~\cite{demmelgrigori2019improved},\cite{MartinssonTroppacta}.
} 
\begin{align*}
\|E_{\rm HMT}\| &= \|(I-\mathcal{P}_{AX})A(I-XM_X^T)\| \leq  \|A(I-XM_X^T)\|. 
\end{align*}
Now let $V\in\mathbb{R}^{n\times \hat r}$ be the leading $\hat r (\leq r)$ right singular vectors of $A$, and choose $M_X=(V^TX)^\dagger V^T$ so that $XM_X^T=X(V^TX)^\dagger V^T=\mathcal{P}_{X,V}$ is an oblique projection onto a subspace of $X$ (the row space is that of $V^T$; unlike other projections in this paper, here $V^TX$ is assumed to have full \emph{row} rank, not column rank). 
We thus have 
$V^T(I-\mathcal{P}_{X,V})=0$, so 
$A(I-\mathcal{P}_{X,V})=A(I-VV^T)(I-\mathcal{P}_{X,V})$. Thus 
\begin{equation}
  \label{eq:Ebound2}
\|E_{\rm HMT}\| = \|A(I-VV^T)(I-\mathcal{P}_{X,V})\|=
\|\Sigma_2(I-\mathcal{P}_{X,V})\|, 
\end{equation}
where $\Sigma_2=\mbox{diag}(\sigma_{\hat r+1},\ldots,\sigma_n)$, so $\|\Sigma_2\| = \|A-A_{\hat r}\|$. 
It turns out that this expression is enough to qualitatively understand why 
the error is near-optimal:  $\|E_{\rm HMT}\|=\|\Sigma_2(I-\mathcal{P}_{X,V})\|
\leq \|\Sigma_2\|\|(I-\mathcal{P}_{X,V})\|_2=\|\Sigma_2\|\|\mathcal{P}_{X,V}\|_2$ by~\eqref{eq:projs}, 
and 
\begin{equation}  \label{eq:PXV}
\|\mathcal{P}_{X,V}\|_2=\|X(V^TX)^\dagger V^T\|_2=\|X(V^TX)^\dagger\|_2
\leq \|X\|_2\|(V^TX)^\dagger\|_2. 
\end{equation}
Now suppose $X$ is Gaussian, hence so is $V^TX\in\mathbb{R}^{\hat r\times r}$. 
Thus by the M-P rule, 
$\|X\|\approx \sqrt{m}+\sqrt{r}$, and 
 $\|(V^TX)^\dagger\|_2\approx 1/(\sqrt{r}-\sqrt{\hat r})$. Thus 
$\|\mathcal{P}_{X,V}\|_2\lesssim \frac{\sqrt{m}+\sqrt{r}}{\sqrt{r}-\sqrt{\hat r}}$. 
It thus follows that with high probability  $  \|E_{\rm HMT}\|\lesssim \frac{\sqrt{m}+\sqrt{r}}{\sqrt{r}-\sqrt{\hat r}}\|A-A_{\hat r}\|$, indicating optimality of HMT up to the factor $\frac{\sqrt{m}+\sqrt{r}}{\sqrt{r}-\sqrt{\hat r}}$ and oversampling $r-\hat r$. 

One might argue that the presence of the large dimension $m$ in this bound is undesirable, and this is indeed an overestimate. 
To get a precise bound some more work is needed, and for this we focus on the Frobenius norm (see~\cite[Cor. 10.10]{halko2011finding} for the 2-norm).  
\begin{theorem}\label{thm:hmt}
Suppose $X,Y$ are Gaussian. Then for any $\hat r\leq r-2$, the HMT error satisfies 
\begin{equation}  \label{eq:HMTbound}
\mathbb{E}\|E_{\rm HMT}\|_F\leq \sqrt{\mathbb{E}\|E_{\rm HMT}\|_F^2}= \sqrt{1+\frac{r}{r-\hat r-1}}\|A-A_{\hat r}\|_F. 
\end{equation}  
\end{theorem}
{\sc proof.}
The first inequality is Cauchy-Schwarz. 
To examine $\mathbb{E}\|E_{\rm HMT}\|_F^2$, we return to~\eqref{eq:Ebound2} and  note that 
$A(I-VV^T)$ and $A(I-VV^T)\mathcal{P}_{X,V}$ lie in complementary row spaces, as the first lies in $V_\perp^T$ and second in $V^T$. It follows that 
\begin{align}
& \|E_{\rm HMT}\|^2_F=\|A(I-VV^T)(I-\mathcal{P}_{X,V})\|_F^2 =
\|A(I-VV^T)\|_F^2+\|A(I-VV^T)\mathcal{P}_{X,V}\|_F^2  \nonumber\\
&=\|\Sigma_2\|_F^2+\|\Sigma_2\mathcal{P}_{X,V}\|_F^2 
= \|\Sigma_2\|_F^2+\|\Sigma_2 (V_\perp^TX)(V^TX)^\dagger V^T\|_F^2.\label{eq:hmtuse}
\end{align}
Now if $X$ is Gaussian then $V_\perp^TX\in\mathbb{R}^{(n-\hat r)\times r}$ and $V^TX\in\mathbb{R}^{\hat r \times r}$ are independent Gaussian. By~\cite[p.~274]{halko2011finding}, it follows that the expected norm of the second term is $\mathbb{E}\|\Sigma_2(V_\perp^TX)(V^TX)^\dagger\|_F^2=\frac{r}{r-\hat r-1}\|\Sigma_2\|_F^2$, so 
\begin{equation}  \label{eq:HMT2}
\mathbb{E}\|E_{\rm HMT}\|_F^2= \left(1+\frac{r}{r-\hat r-1}\right)\|\Sigma_2\|_F^2.
\end{equation}
\qed

The above theorem recovers the bound in~\cite[Thm.~10.5]{halko2011finding} and shows that HMT is optimal to within a (small) oversampling $r-\hat r$ and factor $\sqrt{1+\frac{r}{r-\hat r-1}}$. 

\subsection{Plain generalized \Nystromm}
We now turn to the GN error, $\|E_{\rm GN}\|:=\|A-AX(Y^T\!AX)^{\dagger}Y^T\!A\|$. 
We again start with~\eqref{eq:Eform}, and 
choose  $M_X$ such that $XM_X^T=X(V^TX)^\dagger V^T=\mathcal{P}_{X,V}$
to obtain $\|E_{\rm GN}\| = \|(I-\mathcal{P}_{AX,Y})A(I-\mathcal{P}_{X,V})\|$. 
A difference from HMT is that $\mathcal{P}_{AX,Y}$ is not an orthogonal projector. We can still bound $\|E_{\rm GN}\|$ as 
\begin{align}
\|E_{\rm GN}\| &= \|(I-\mathcal{P}_{AX,Y})A(I-\mathcal{P}_{X,V})\| = \|(I-\mathcal{P}_{AX,Y})A(I-VV^T)(I-\mathcal{P}_{X,V})\| \label{eq:EGN}\\
& \leq \|(I-\mathcal{P}_{AX,Y})\|_2\|A(I-VV^T)(I-\mathcal{P}_{X,V})\|
=\|\mathcal{P}_{AX,Y}\|_2  \|E_{\rm HMT}\|,\nonumber
\end{align}
where $E_{\rm HMT}$ is the HMT error as in~\eqref{eq:Ebound2}. 
Hence generalized \Nystrom has error within 
$\|\mathcal{P}_{AX,Y}\|_2$ of HMT. 
When $X,Y$ are taken to be independent Gaussian matrices we have 
$\|\mathcal{P}_{AX,Y}\|_2\lesssim \frac{\sqrt{n}+\sqrt{r+\ell}}{\sqrt{r+\ell}-\sqrt{r}}$; this can be shown by the same argument after~\eqref{eq:PXV}, by noting that $\mathcal{P}_{AX,Y}=\mathcal{P}_{Q,Y}$ where $QR=AX$ is the QR factorization. The near-optimality of GN follows. 

Again, with more work one can get a precise bound that does not involve $\sqrt{m},\sqrt{n}$. The following reproduces~\cite[Thm.~4.3]{tropp2017practical}. 
\begin{theorem}\label{thm:plainGN}
Suppose $X,Y$ are Gaussian. 
Then for any $\hat r\leq r-2$, the error with plain generalized \Nystrom\eqref{eq:start} satisfies 
\begin{equation}  \label{eq:gnbound}
\mathbb{E}\|E_{\rm GN}\|_F\leq \sqrt{\mathbb{E}\|E_{\rm GN}\|_F^2}\leq \sqrt{1+\frac{r+\ell}{\ell-1}}\sqrt{1+\frac{r}{r-\hat r-1}}\|A-A_{\hat r}\|_F
=\sqrt{1+\frac{r+\ell}{\ell-1}}\sqrt{\mathbb{E}\|E_{\rm HMT}\|_F^2}.
\end{equation}  
\end{theorem}
{\sc proof.}
Following~\cite[\S A.3]{tropp2017practical} we write 
$\mathcal{P}_{AX,Y}A=Q(Q^T+Z )A$, 
where $Q=\mbox{orth}(AX)$, 
so that 
$E_{\rm GN}=(I-\mathcal{P}_{AX,Y})A=(I-QQ^T)A+QZ A$. Noting that $(I-QQ^T)A=Q_\perp Q_\perp^TA$ is precisely the error from HMT, and
$QZ A$ can be expressed as
\[
QZ A = Q((Y^TQ)^\dagger Y^T-Q^T)A
 = Q((Y^TQ)^\dagger Y^T-Q^T)Q_\perp Q_\perp^TA
=Q(Y^TQ)^\dagger (Y^TQ_\perp) Q_\perp^TA
\]
where the second equality holds because 
$((Y^TQ)^\dagger Y^T-Q^T)Q=0$. Since $\|Q^T_\perp A\|=\|E_{\rm HMT}\|$, and if $Y$ is Gaussian then $Y^TQ$ and  $Y^TQ_\perp$ are independent Gaussian matrices, so using~\cite[p.~274]{halko2011finding} again with Pythagoras, 
we obtain~\eqref{eq:gnbound}.
\qed

In words, the expected GN error is optimal to within the factor $\sqrt{\left(1+\frac{r}{r-\hat r-1}\right)\left(1+\frac{r+\ell}{\ell-1}\right)}$, and at most $\sqrt{1+\frac{r+\ell}{\ell-1}}$ times worse than that of HMT. 


\subsection{Stabilized GN}
We turn to stabilized GN~\eqref{eq:ournystrom} and derive bounds on the error $\|E_{\rm SGN}\|:=\|A- AX(Y^T\tilde AX)_\epsilon^{\dagger}Y^T\!A\|$. Roughly, the goal is to show $\|E_{\rm SGN}\|$ is on the same order as 
 with plain GN. We will show this is true for any $\tilde A=A+O(u\|A\|)$.
Our analysis below is unfortunately rather lengthy. 
Nonetheless, the results will be needed for the stability analysis of SGN, and also to highlight the potential instability of plain GN. 

In SGN, instead of $\P_{AX,Y}A$ as in GN, the approximant is 
\begin{align}
AX(Y^T\tilde AX)_\epsilon^\dagger Y^T\!A
=(AX(Y^T\tilde AX)_\epsilon^\dagger Y^T)A=:\tilde \P_{AX,Y}A. \label{eq:tildePAXY}
\end{align}
Note that 
$\tilde \P_{AX,Y}$ is not necessarily a projector since $\tilde\P_{AX,Y}^2\neq \tilde\P_{AX,Y}$, 
and unlike $\P_{AX,Y}$, $\tilde\P_{AX,Y}$ does not project onto the range of $AX$, because $\tilde\P_{AX,Y}AX\neq AX$. 
This is why the arguments for GN above do not carry over directly. 

Here and below, we use $\O(1)$ to suppress terms like $\sqrt{m},\sqrt{n},r$ (but of course not $1/\epsilon,1/\sigma_r(A)$, etc). This might seem 
like an oversimplification, but this is standard practice in stability analysis (e.g.~\cite{nahi11})
and the hidden constants below are attached to $\epsilon$ (in the order of machine precision)
rather than (the usually much larger) $\|A-A_r\|$. 
The $O(1)$ notation does not hide such terms and continues to keep them separate. We say that $\|X\|=O(1)$ holds with exponentially high probability if 
$\mathbb{P}[\|X\|\geq t]\leq \exp(-ct)$ for some $c>0$. 

\ignore{
\hrule 
We note that because $\epsilon=O(u)$, there exists $\At: = A+\Delta A$ with $\|\Delta A\|=O(u)$ such that 
\[AX(Y^T\!AX)_\epsilon^\dagger Y^T\!A = 
\At X(Y^T\At X)^\dagger Y^T\At, \]
that is, the approximant with $\epsilon$-truncation is the exact approximant 
of a perturbed matrix $\At$ without the $\epsilon$-truncation. 

We now invoke the results from the previous section to obtain 
$
\At X(Y^T\At X)^\dagger Y^T\At = \At + E, 
$ where the error $E$ has been worked out to be $O(\|\widetilde\Sigma_2\|)$. Since singular values are well-conditioned, we have $\|\widetilde\Sigma_2\|=\|\Sigma_2\|+\|\Delta A\|=\|\Sigma_2\|+O(u)$. 
It follows that the $\epsilon$-truncation is a safe operation. 
We believe the analysis here is significantly simpler than \cite{chiu2013sublinear}. 

(here we do not keep track of constants that depend on $m,n,r$; we write $O(u)$ to mean a term that is $u$ times a modest polynomial of $m,n,r$. )
\hrule 
}



\subsubsection{Lemmas}
We start with three lemmas. 
\begin{lemma}\label{lem:g2bound}
  Let $G$ be $m\times n$ Gaussian with $m-1\geq n\geq 2$. Then 
  \begin{equation}
    \label{eq:Gau}
  \mathbb{E}  \|G^\dagger\|_2^2\leq 
\frac{e^2m}{(m-n)^2-1}.
  \end{equation}
\end{lemma}
{\sc proof.}
The argument closely follows \cite[Prop. A.4]{halko2011finding}, which bounds the first moment of the pseudoinverse of $m\times n$ Gaussians 
$\mathbb{E}  \|G^\dagger\|<\frac{e\sqrt{m}}{m-n}$.
We use the inequality below, which follows from~\cite{chen2005condition}, \cite[Prop. A.3]{halko2011finding}:
\begin{equation}  \label{eq:Gt}
\mathbb{P}[\|G^\dagger\|_2^2>t]\leq \frac{1}{\sqrt{2\pi(m-n+1)}}\left(\frac{e\sqrt{m}}{m-n+1}\right)^{m-n+1}t^{-(m-n+1)/2}. 
\end{equation}
Now write $C:=\frac{1}{\sqrt{2\pi(m-n+1)}}\left(\frac{e\sqrt{m}}{m-n+1}\right)^{m-n+1}$. For any $E>0$ we have 
\begin{align*}
\mathbb{E}\|G^\dagger\|_2^2
&= \int_0^\infty \mathbb{P}[\|G^\dagger\|_2^2>t]dt
\leq E + \int_E^\infty \mathbb{P}[\|G^\dagger\|_2^2>t]dt\\
&\leq E + C\int_E^\infty t^{-(m-n+1)/2}dt=E +\frac{1}{(m-n+1)/2-1} C E^{-(m-n+1)/2+1}. 
\end{align*}
We minimize this expression with respect to $E$, which gives 
$1=CE^{-(m-n+1)/2}$. 
With this choice, we obtain 
\begin{align*}
\mathbb{E}\|G^\dagger\|_2^2 &\leq 
E(1+\frac{1}{(m-n-1)/2})
=C^{\frac{2}{m-n+1}}(1+\frac{1}{(m-n-1)/2})\\
&=\left(\frac{1}{2\pi(m-n+1)}\right)^{\frac{1}{m-n-1}}\left(\frac{e\sqrt{m}}{m-n+1}\right)^2(1+\frac{2}{m-n-1})\\
&\leq \frac{e^2m}{(m-n)^2-1}.
\end{align*}
\qed 

In view of the Marchenko-Pastur rule we expect  $\|G^\dagger\|_2^2=1/(\sigma_{\min}(G))^2\approx 1/(\sqrt{m}-\sqrt{n})^2$, so~\eqref{eq:Gau} is a reasonable bound, though undoubtedly the constant can be improved.

\begin{lemma}\label{lem:AXX}
Let $A,X,Y$ be such that $X,Y$ are Gaussian, $AX$ is full column rank and $Y^T\!AX$ is tall, and let $\At$ be any matrix such that $\At=A+\delta A$ where $\|\delta A\|=\O(u\|A\|)$. 
Suppose also that $r\geq 2$ and $\ell\geq 1$.  
 Then with exponentially high probability 
\begin{equation}  \label{eq:AX1}
\|A X(Y^T\!A X)^\dagger\|=O(1),\quad 
\|A X(Y^T\At X)_\epsilon^\dagger\|=O(1). 
\end{equation}  
Suppose further that $\epsilon$ is chosen s.t. $\|Y^T\!\delta AX\|_2\leq \epsilon$. Then 
\begin{equation}  \label{eq:specalways}
\|Y^T\!A X(Y^T\!A X+\varepsilon)_\epsilon^\dagger \|_2
\leq 2  , 
\end{equation}
(which holds deterministically), and 
\begin{equation}  \label{eq:AX13}
\mathbb{E}\|A X(Y^T\!A X)^\dagger\|_2^2\leq 
\frac{e^2(r+\ell)}{\ell^2},
\quad 
\mathbb{E}\|A X(Y^T\At X)_\epsilon^\dagger\|_2^2\leq 
\frac{4e^2(r+\ell)}{\ell^2}. 
\end{equation}  
\end{lemma}
{\sc proof.} 
For the first statement, let $AX = U\Sigma V^T$ be the SVD. 
\[A X(Y^T\!A X)^\dagger =U\Sigma V^T(Y^TU\Sigma V^T)^\dagger 
=U\Sigma V^TV(Y^TU\Sigma )^\dagger,
=U\Sigma (Y_1\Sigma )^\dagger, \]
where $Y_1:=Y^TU$.
Hence it suffices to prove $\|\Sigma (Y_1\Sigma )^\dagger\|=O(1)$. 
For this we see that $\|Y_1\Sigma (Y_1\Sigma )^\dagger\|= 1$, and 
$\|\Sigma (Y_1\Sigma )^\dagger\|=
\|Y_1^\dagger(Y_1\Sigma (Y_1\Sigma )^\dagger)\|
\leq \|Y_1^\dagger\|=O(1)$, which follows from the fact that $Y_1$ is tall-Gaussian, hence well-conditioned (by M-P $\|Y_1^\dagger\|\approx \frac{1}{\sqrt{r+\ell}-\sqrt{r}}$).  Furthermore, 
from $\|A X(Y^T\!A X)^\dagger\| = \|Y_1^\dagger\| $
and Lemma~\ref{lem:g2bound} 
we obtain the first statement in~\eqref{eq:AX13}:
\begin{equation}  \label{eq:Y1dagger2}
\mathbb{E}\|A X(Y^T\!A X)^\dagger\|_2^2 = \mathbb{E}\|Y_1^\dagger\|_2^2
\leq \frac{e^2(r+\ell)}{\ell^2}. 
\end{equation}


For the second claim in~\eqref{eq:AX1}, 
first note that the above proof is inapplicable because for any fixed $X,Y$, there exists a small perturbation $\tilde A$ of $A$ such that $Y_1$ is ill-conditioned. 
Instead, we write 
$A X(Y^T\At X)_\epsilon^\dagger = A X(Y^T\!A X+\varepsilon)_\epsilon^\dagger $
where $\|\varepsilon\|_2\leq \epsilon$ by assumption, 
and note that 
$\|(Y^T\!A X+\varepsilon) (Y^T\!A X+\varepsilon)_\epsilon^\dagger \|_2\leq 1$.
Thus by the triangle inequality 
\[
\|Y^T\!A X(Y^T\!A X+\varepsilon)_\epsilon^\dagger \|
=
\|(Y^T\!A X+\varepsilon-\varepsilon)(Y^T\!A X+\varepsilon)_\epsilon^\dagger \|
\leq 
O(1)+\|\varepsilon (Y^T\!A X+\varepsilon)_\epsilon^\dagger \|=  O(1),
\]
because $\|(Y^T\!A X+\varepsilon)_\epsilon^\dagger\|_2\leq 1/\epsilon$. 
If further $\|Y^T\!\delta AX\|_2\leq \epsilon$, then in the spectral norm in the final inequality both terms are bounded deterministically by 2, giving~\eqref{eq:specalways}. 

Again using the SVD $AX=U\Sigma V^T$, 
\[
Y^T\!A X(Y^T\!A X+\varepsilon)_\epsilon^\dagger
=(Y^TU)\Sigma V^T(Y^T\!A X+\varepsilon)_\epsilon^\dagger, 
\]
so 
\begin{align}
\|A X(Y^T\!A X+\varepsilon)_\epsilon^\dagger\|&=
\|\Sigma V^T(Y^T\!A X+\varepsilon)_\epsilon^\dagger\|\nonumber\\
&=\|(Y^TU)^\dagger (Y^T\!A X(Y^T\!A X+\varepsilon)_\epsilon^\dagger) \| \nonumber\\
&\leq\|Y_1^\dagger\| \|(Y^T\!A X(Y^T\!A X+\varepsilon)_\epsilon^\dagger) \| =O(1),\label{eq:AXYnorm}
\end{align}
because again $Y_1$ is a tall Gaussian matrix hence well-conditioned with high probability. 
Finally, we bound the expected value in the spectral norm, using~\eqref{eq:specalways}, as 
\begin{align*}
\mathbb{E}\|A X(Y^T\!A X+\varepsilon)_\epsilon^\dagger\|_2^2
&\leq \mathbb{E}\left(\|Y_1^\dagger\|_2^2 \|(Y^T\!A X(Y^T\!A X+\varepsilon)_\epsilon^\dagger) \|_2^2\right)
\leq \mathbb{E}(4\|Y_1^\dagger\|_2 ^2)
\leq 
\frac{4e^2(r+\ell)}{\ell^2},
\end{align*}
completing the proof. \qed 
\ignore{
The final (third) statement in~\eqref{eq:AX1} is now immediate: 
\begin{align*}
\|\At X(Y^T\At X)_\epsilon^\dagger\|&\leq \|A X(Y^T\At X)_\epsilon^\dagger\|+\|\epsilon X(Y^T\At X)_\epsilon^\dagger\|\\
&\leq O(1)+\|\epsilon\| \|X\|\|(Y^T\At X)_\epsilon^\dagger\|=O(1),  
\end{align*}
because $\|(Y^T\At X)_\epsilon^\dagger\|\leq 1/\epsilon$ by definition. 
}




Let us note that it is not always true that 
$\|A X(Y^T\At X)^\dagger\|=O(1)$, because 
assuming $Y^T\!AX$ has $O(u)$ singular values (which is typically the case), there exist perturbations $\delta A$ such that 
$Y^T(A+\delta A)X$ has singular values $\ll \epsilon$. 
This fact will be important in Section~\ref{sec:gennyst}.
We also note that we have not proved $\|(Y^T\!A X)^\dagger Y^T\!A\|=O(1)$, which appears to be nontrivial. 

Below and in Section~\ref{sec:gennyst}, 
we use $\eepsilon$ to denote either a matrix or a scalar such that $\|\eepsilon\|=\O(u)$. 
The precise value of $\eepsilon$ may change from appearance to appearance; this follows the practice in e.g.~\cite{nahi11}. This notation simplifies the analysis without losing the essense. 


We now turn to a key result that allow us to reduce the analysis of SGN to one for plain GN. 
\begin{lemma}
With $\tilde\P_{AX,Y}=AX(Y^T\tilde AX)_\epsilon^\dagger Y^T$ as in \eqref{eq:tildePAXY}, 
\begin{equation}  \label{eq:AXp}
\tilde\P_{AX,Y}AX=AX+\eepsilon . 
\end{equation}  
\end{lemma}
{\sc proof.}
First note that because $\epsilon=\O(u)$, there exists $\At_\epsilon: 
= \At+\eepsilon$, 
which can also be written as $A+\Delta A$ with $\|\Delta A\|=\O(u)$, such that 
\[
(Y^T\!\At X)_\epsilon^\dagger = (Y^T\At_\epsilon X)^\dagger, 
\]
that is, the $\epsilon$-pseudoinverse is the exact pseudoinverse of the perturbed matrix $\At_\epsilon$. 
Hence we have $\tilde\P_{AX,Y}AX=AX(Y^T\!\At X)_\epsilon^\dagger (Y^T\!AX)=AX(Y^T\At_\epsilon X)^\dagger (Y^T\!AX)$. Now 
\begin{align*}
\tilde\P_{AX,Y}AX&=AX(Y^T\At_\epsilon X)^\dagger (Y^T\!AX)
=(\At_\epsilon+\eepsilon)X(Y^T\At_\epsilon X)^\dagger (Y^T(\At_\epsilon+\eepsilon)X)\\
&=\At_\epsilon X(Y^T\At_\epsilon X)^\dagger (Y^T\At_\epsilon X)+\eepsilon
=\P_{\At_\epsilon X}\At_\epsilon X+\eepsilon,
  \end{align*}
where for the penultimate equality we used the fact
$\|X(Y^T\At_\epsilon X)_\epsilon^\dagger Y^T\At_\epsilon X\|\leq 
\|X\|\|(Y^T\At_\epsilon X)_\epsilon^\dagger Y^T\At_\epsilon X\|_2=\|X\|=O(1)$, and 
$\|\At_\epsilon X(Y^T\At_\epsilon X)_\epsilon^\dagger\|=O(1)$ from Lemma~\ref{lem:AXX}. 
Letting $Y^T\At_\epsilon X=\Ut\St \Vt^T
=[\Ut_1,\Ut_2]\big[
\begin{smallmatrix}
\St_1 & \\   & \St_2
\end{smallmatrix}
\big]
\big[
\begin{smallmatrix}
\Vt_1^T\\\Vt_2^T
\end{smallmatrix}
\big]$ be the SVD where $\|\St_2\|\leq \epsilon$, we have 
\begin{align*}
\P_{\At_\epsilon X}\At_\epsilon X=\At_\epsilon X
[\Vt_1, \Vt_2]
  \begin{bmatrix}
    I & \\ & 0 
  \end{bmatrix}
[\Vt_1, \Vt_2]^T
=\At_\epsilon X\Vt_1\Vt_1^T=\At_\epsilon X-\At_\epsilon X\Vt_2\Vt_2^T. 
\end{align*}
Since $\At_\epsilon X=AX+\eepsilon$, it suffices to show that $\At_\epsilon X\Vt_2=\eepsilon$, or equivalently that $A X\Vt_2=\eepsilon$. 

Let $A X=
\Uh\Sh\Vh^T
$ 
be the SVD. Then 
$Y^TAX=\Omega\Sh\Vh^T$ where $\Omega$ is tall Gaussian,  hence well-conditioned. 
$Y^TAX \Vt_2 = \Omega\Sh\Vh^T\Vt_2 $, and so 
$(\eepsilon =) \|Y^TAX \Vt_2\| = \|\Omega\Sh\Vh^T\Vt_2\|\leq \|\Omega^\dagger\| \|\Sh\Vh^T\Vt_2\| = O(\|\Sh\Vh^T\Vt_2\|)$. It follows that 
$\|\Sh\Vh^T\Vt_2\|=\eepsilon$, and hence 
$\|AX\Vt_2\| = \|\Sh\Vh^T\Vt_2\|=\eepsilon$, as required. 
\qed




We note that in the above lemma and~\eqref{eq:AX1}, the Gaussianity of $X,Y$ is not essential; any class of random matrix such that the entries are $O(1)$ and a rectangular realization is well-conditioned would suffice, including the SRFT and SRHT matrices. 


\subsubsection{Accuracy of stabilized generalized \Nystromm}
Now we turn to the main subject of 
assessing the accuracy of stabilized generalized \Nystromm. 
As before, our goal is to bound $\|E_{\rm SGN}\|:=\|A-\hat A_r\|=\|A-\tilde\P_{AX,Y}A\|=\|(I-\tilde\P_{AX,Y})A\|$. 

From~\eqref{eq:AXp} we have $(I-\tilde\P_{AX,Y})AX=\eepsilon$, and so
$(I-\tilde\P_{AX,Y})A=(I-\tilde\P_{AX,Y})A(I-XM_X)+\eepsilon$ for any $M_X$ such that $\|M_X\|=O(1)$. We take $M_X=(V^TX)^\dagger V^T$ as before so that $XM_X=\mathcal{P}_{X,V}$, which satisfies $\|M_X\|=O(1)$. 
This yields
\begin{align*}
\|E_{\rm SGN}\|=\|(I-\tilde\P_{AX,Y})A\|&=\|(I-\tilde\P_{AX,Y})A(I-\mathcal{P}_{X,V})\|+\eepsilon, 
\end{align*}
where the first term can be bounded as in~\eqref{eq:EGN} and~\eqref{eq:gnbound}:
\begin{align}
\|(I-\tilde\P_{AX,Y})A(I-\mathcal{P}_{X,V})\|
&=\|(I-\tilde\P_{AX,Y})A(I-VV^T)(I-\mathcal{P}_{X,V})\|\nonumber\\
&\leq \|I-\tilde \P_{AX,Y}\|_2 \| E_{\rm HMT}\|.\label{eq:sgnuse}
\end{align}
We have used the inequality $\|AB\|\leq \|A\|_2\|B\|$, which holds for any unitarily invariant norm~\cite[Cor.~3.5.10]{hornjohntopics}. 
Note that we cannot take $\|I-\tilde \P_{AX,Y}\|_2=\|\tilde \P_{AX,Y}\|$, as $\tilde \P_{AX,Y}$ is not exactly a projection (it is an approximate projection $\tilde \P_{AX,Y}^2=\tilde \P_{AX,Y}+\eepsilon$). We can still bound 
$\|I-\tilde \P_{AX,Y}\|_2 \leq 1+\|\tilde \P_{AX,Y}\|_2$, and 
\begin{align*}
\|\tilde \P_{AX,Y}\|&=\|AX(Y^T\tilde AX)_\epsilon^\dagger Y^T\|
\leq \|AX(Y^T\tilde AX)_\epsilon^\dagger\| \|Y^T\|=O(1), 
\end{align*}
where we used $\|AX(Y^T\tilde AX)_\epsilon^\dagger\|=O(1)$ from Lemma~\eqref{lem:AXX}. Together with~\eqref{eq:sgnuse} we conclude that 
$E_{\rm SGN} = O(1)E_{\rm HMT}$. Note, however, that this $O(1)$ notation suppresses terms like $\sqrt{m}$. 
Once again, a more precise bound can be obtained, as follows. 

\begin{theorem}\label{thm:sgnmain}
Let the assumptions in Lemma~\ref{lem:AXX} be satisfied, including $\|Y^T\!\delta AX\|_2\leq \epsilon$, where $\delta A = \At-A$.  
For any $\hat r\leq r-2$, the error with stabilized generalized \Nystrom\eqref{eq:start} for any $\tilde A$ satisfies 
  \begin{equation}\label{eq:sgnmain}
    \begin{split}
\mathbb{E}\|E_{\rm SGN}\|_F    &\leq  
\frac{2\sqrt{e}(r+\ell)}{\ell}\sqrt{1+\frac{r}{r-\hat r-1}}\|A-A_{\hat r}\|_F+\eepsilon\\
&=\frac{2\sqrt{e}(r+\ell)}{\ell}\sqrt{\mathbb{E}\|E_{\rm HMT}\|_F^2   }+\eepsilon.
    \end{split}
  \end{equation}
\end{theorem}

{\sc proof}. Denote by $U_\perp,V_\perp$ the trailing singular vectors of $A$, from the $(\hat r+1)$th. We have 
\begin{align}
\|E_{\rm SGN}\|    =& \|(I-\tilde\P_{AX,Y})A(I-\mathcal{P}_{X,V})\|+\eepsilon
=\|(I-\tilde\P_{AX,Y})U_\perp U_\perp^TAV_\perp V_\perp^T(I-\mathcal{P}_{X,V})\|+\eepsilon\nonumber\\
\leq &\|U_\perp U_\perp^TAV_\perp V_\perp^T(I-\mathcal{P}_{X,V})\|
+\|\tilde\P_{AX,Y}U_\perp U_\perp^TAV_\perp V_\perp^T(I-\mathcal{P}_{X,V})\|+\eepsilon
\nonumber\\
=&:\|E_{\rm HMT}\|+\|E_1\|+\eepsilon\label{eq:sgnfirst}
\end{align}
where 
\begin{align}
\|E_1\|:&=\|\tilde\P_{AX,Y}U_\perp U_\perp^TAV_\perp V_\perp^T(I-\mathcal{P}_{X,V})\|
=\|AX(Y^T\!\At X)_\epsilon^\dagger Y^T U_\perp U_\perp^TAV_\perp V_\perp^T(I-\mathcal{P}_{X,V})\|\nonumber
\\
&\leq \|AX(Y^T\!\At X)_\epsilon^\dagger \|_2\|Y^T U_\perp U_\perp^TAV_\perp V_\perp^T(I-\mathcal{P}_{X,V})\|.
  \label{eq:E1}
\end{align}
Hence by Cauchy-Schwarz 
and Lemma~\ref{lem:AXX} we obtain 
\begin{align}
\mathbb{E}\|E_1\|_F^2&\leq
\mathbb{E}\|AX(Y^T\!\At X)_\epsilon^\dagger \|_2^2
\ \mathbb{E}\|Y^T U_\perp U_\perp^TAV_\perp V_\perp^T(I-\mathcal{P}_{X,V})\|_F^2\nonumber 
\\&\leq 
\frac{4e^2(r+\ell)}{\ell^2}
\mathbb{E}\|Y^T U_\perp U_\perp^TAV_\perp V_\perp^T(I-\mathcal{P}_{X,V})\|_F^2. \label{eq:firstterm}
  \end{align}
For the second term, we take the expectations with respect to $X$ and $Y$ separately. Namely 
$M_X:=U_\perp U_\perp^TAV_\perp V_\perp^T(I-\mathcal{P}_{X,V})$ is independent of $Y$, and for any fixed $M_X$, by~\cite[Prop. 10.1]{halko2011finding} we have $\mathbb{E}\|YM_X\|_F^2 = \sqrt{r+\ell}\|M_X\|_F^2$. Therefore, now making the random variables explicit and noting that $\|M_X\|$ is equal to the HMT error, 
\begin{align}
\mathbb{E}&\|Y^T U_\perp U_\perp^TAV_\perp V_\perp^T(I-\mathcal{P}_{X,V})\|_F^2
=\mathbb{E}_{X,Y}\|Y^T M_X\|_F^2 =(r+\ell)(\mathbb{E}_{X}\|M_X\|_F)^2 \nonumber \\
&= (r+\ell)(\mathbb{E}_{X}\|E_{\rm HMT}\|_F)^2 \leq \|\Sigma_2\|_F^2(r+\ell)(1+\frac{r}{r-\hat r-1}), \label{eq:HMTF}
\end{align}
where we used~\eqref{eq:hmtuse} for the final inequality. 
Substitute this into~\eqref{eq:sgnfirst} together with \eqref{eq:E1} and \eqref{eq:firstterm} to obtain the required bound. \qed

As the analysis shows, and as mentioned after~\eqref{eq:ournystrom}, 
$(Y^T\!\At X)_\epsilon^\dagger$ can be the pseudoinverse of 
any matrix that is $\epsilon$-close to $Y^T\!AX$ for the results to hold.


\section{Numerical stability of generalized \Nystromm}\label{sec:gennyst}
We now examine the numerical stability of plain and stabilized GN, taking into account roundoff errors in floating-point arithmetic. We will first establish the stability of stabilized GN, then discuss the potential instability of plain GN.

In the remainder of this section, without loss of generality we assume $\|A\|_2=1$, so that any ill-conditioning or approximability by low-rank matrices comes from the presence of small singular values rather than large ones. This is not a fundamental assumption but helps simplify the arguments. 
We also (continue to) assume that $X,Y$ are rectangular Gaussian random matrices so that $\|X\|,\|Y\|,\|X^\dagger\|,\|Y^\dagger\|$ are all $\O(1)$. 
\subsection{Stability of stabilized GN}\label{sec:staGN}
Consider the stabilized GN $AX(Y^T\!AX)_\epsilon^\dagger Y^T\!A$ in finite-precision arithmetic. Recall that this is what we attempt to compute rather than 
$AX(Y^T\!\At X)_\epsilon^\dagger Y^T\!A$ in~\eqref{eq:ournystrom}; the latter was analyzed above for theoretical use below. 
We assume that $AX(Y^T\!AX)_\epsilon^\dagger$ is computed first, then $Y^T\!A$ is multiplied. 
Note that each row of $AX(Y^T\!AX)_\epsilon^\dagger$ involves solving an underdetermined system of linear equations with respect to $\epsilon$-pseudoinverses. 

Below we will use~\cite[Thm 21.4]{Higham:2002:ASNA} (or~\cite{demmel1993improved}), which states that a numerically full-rank underdetermined linear system 
$x=A^\dagger b$ solved by the QR-based (or SVD-based) method is backward stable, i.e., the computed solution $\hat x$ is the minimum-norm solution for a perturbed problem $\hat x=(A+\delta A)^\dagger b$, where $\|\delta A\|_2
=\eepsilon \\A\|_2$. 

We now state the main stability result of SGN. 
Let us denote the $i$th row of a matrix $Z$ by $[Z]_i$. 
 \begin{theorem}\label{thm:sta}
Suppose that 
$AX(Y^T\!AX)_\epsilon^\dagger Y^TA$ is computed as 
$AX(Y^T\!AX)_\epsilon^\dagger$ times $Y^TA$, and each row of $AX(Y^T\!AX)_\epsilon^\dagger$ is computed by a backward stable underdetermined linear solver. 
Then for some $\At_i=A+\delta A_i$ with $\|\delta A_i\|\leq \eepsilon\|A\|$, for every $i$ 
\begin{equation}
  \label{eq:ff}
\|[fl(AX(Y^T\!AX)_\epsilon^\dagger Y^TA) - A]_i\|_2= 
\|[(AX(Y^T\!\tilde A_iX)_\epsilon^\dagger Y^TA) - A]_i\|_2+\eepsilon
=\|[E_{\rm SGN}]_i\|_2+\eepsilon, 
\end{equation}
Suppose further that $X,Y$ are Gaussian and $\epsilon \geq \|Y^T\!(A-\At_i) X\|_2$ for all $i$ (but still $\epsilon=\O(u\|A\|)$). Then 
\begin{equation}  \label{eq:Efloat}
\mathbb{E}\|fl(AX(Y^T\!AX)_\epsilon^\dagger Y^TA) - A\|_F
\leq \left(\frac{4e\sqrt{r}(r+\ell)}{\ell}+1\right) \mathbb{E}\|E_{\rm HMT}\|_F+\eepsilon.
\end{equation}
\ignore{
Then for every $i$, the computed approximant of the $i$th row $[AX(Y^T\!AX)_\epsilon^\dagger Y^TA]_i$ satisfies 
\begin{equation}
  \label{eq:fl}
\|[fl(AX(Y^T\!AX)_\epsilon^\dagger Y^TA) - A]_i\|_2= 
\|[(AX(Y^T\!\tilde AX)_\epsilon^\dagger Y^TA) - A]_i\|_2+\epsilon
\end{equation}
for some $\tilde A$ such that $\tilde A=A+\delta A$, $\|\delta A\|=\O(u\|A\|)$. 
}
 \end{theorem}
{\sc proof.}
Define 
$s_i^T=[AX(Y^T\!AX)_\epsilon^\dagger]_i$, 
and let its computed approximant be 
$\hat s_i^T=fl(AX(Y^T\!AX)_\epsilon^\dagger)$. 
Then for each $i$, $s_i$ is the minimum-norm solution to the underdetermined linear system with $\epsilon$-truncated singular values
\[s_i = ((Y^T\!AX)^T)_\epsilon^\dagger [AX]_i^T
=(X^T\!A^TY)_\epsilon^\dagger [AX]_i^T=:M_\epsilon^\dagger [AX]_i^T. \]
It satisfies 
$ M_\epsilon  s_i = [AX]_i^T-\tilde\epsilon$,
where $\tilde\epsilon$ is the component in $(AX)_i^T$ that does not belong to the column space of $M_\epsilon=(X^T\!A^TY)_\epsilon$, which is $\O(u)=\eepsilon$ from~\eqref{eq:AXp}.


We would now invoke \cite[Thm 21.4]{Higham:2002:ASNA}, but 
an issue is that the theorem requires the underdetermined system to be of numerically full row rank, so it is not immediately applicable here. 
Let $s_i = M_\epsilon^\dagger [AX]_i^T$
be computed via a backward stable method to yield $\hat s_i$. 
Then the first step is to project the vector $[AX]_i^T$ onto 
$\hat U$, the computed column space of $M_\epsilon$, which 
is the exact column space of $M+\eepsilon$. 
Now 
defining $\hat b=fl(\hat U^T[AX]_i^T)=\hat U^T[AX]_i^T+\epsilon$, 
the task becomes to solve the 
smaller underdetermined system  $(\hat U^TM)^\dagger  \hat b,  $
 where $\hat U^TM$ is numerically full-rank with singular values $>\epsilon$. Thus by \cite[Thm 21.4]{Higham:2002:ASNA}, its computed solution $\hat s_i$ 
satisfies 
$\hat s_i=(\hat U^TM+\eepsilon)^\dagger  \hat b=(\hat U^TM+\eepsilon)_\epsilon^\dagger  \hat b$, 
from which we obtain
\begin{align}
\hat s_i&
=(\hat U^TM+\eepsilon)^\dagger (\hat U^T[AX]_i^T+\eepsilon) 
 = (M+\epsilon_i)_\epsilon^\dagger ([AX]_i^T+\eepsilon)_\epsilon.      \label{eq:mihatpert}
\end{align}
The last equality holds because
$(\hat U^TM+\eepsilon)^\dagger \hat U^T=
(\hat U\hat U^TM+\eepsilon)^\dagger =:
(M+\epsilon_i)_\epsilon^\dagger$, since 
$\hat U\hat U^TM=M+\eepsilon$ from the definition of $\hat U$. 

Therefore 
for each $i$, we can write 
$[fl(AX(Y^T\!AX)_\epsilon^\dagger)]_i=\hat s_i^T=
[AX+\eepsilon]_i(Y^T\!AX+\epsilon_i)_\epsilon^\dagger
$
and so 
\begin{align}
[fl(AX(Y^T\!AX)_\epsilon^\dagger Y^T\!A)]_i
&=fl([AX+\eepsilon]_i(Y^T\!AX+\epsilon_i)_\epsilon^\dagger Y^T\!A)\nonumber \\
&=[AX]_i(Y^T\!AX+\epsilon_i)_\epsilon^\dagger Y^T\!A
+\eepsilon\|[AX]_i(Y^T\!AX+\epsilon_i)_\epsilon^\dagger\| \|Y^T\!A\|  \nonumber \\
&=[AX]_i(Y^T\!AX+\epsilon_i)_\epsilon^\dagger Y^T\!A
+\eepsilon,\label{eq:break}
\end{align}
where we used~\eqref{eq:AX1} for the last equality, namely
$\|[AX]_i(Y^T\!AX+\epsilon_i)_\epsilon^\dagger Y^T\!A\| \leq \|[AX]_i(Y^T\!AX+\epsilon_i)_\epsilon^\dagger \| \|Y^T\!A\| =\O(1)$.


\ignore{
 From before we also have for the whole matrix 
\[(AX)(Y^T\!AX+\epsilon_i)_\epsilon^\dagger Y^T\!A =
A+E_{\rm SGN},\]
so clearly the every row 
$[AX(Y^T\!AX+\epsilon_i)_\epsilon^\dagger Y^T\!A]_i$ is within 
$[E_{\rm SGN}]_i$
of $A_i$ and by~\eqref{eq:ff} so is 
the computed $fl([AX(Y^T\!AX+\epsilon_i)_\epsilon^\dagger Y^T\!A]_i)$. This shows that for each row, the error with a computed SGN is 
comparable to that of an exact SGN. This statement can be used repeatedly to conclude that the whole matrix $AX(Y^T\!AX)_\epsilon^\dagger Y^T\!A$ is computed with error comparable to that of an exact SGN, that is, 
\begin{equation}
  \label{eq:flbad}
\|fl(AX(Y^T\!AX)_\epsilon^\dagger Y^TA) - A\|_F=
\mathcal{O}(
\|E_{\rm SGN}\|_2
)+\eepsilon.
\end{equation}
However, in the worst case this would incur an error growth by a factor $\sqrt{m}$ (which is hidden by $\mathcal{O}$ here), since the $\epsilon_i$ in~\eqref{eq:break} depends on $i$. This is due to the worst-case nature of the derivation, which appears necessary because $\epsilon_i$ are due to roundoff errors, which cannot necessarily be modelled by random matrices. This $\sqrt{m}$-amplification is unlikely and not observed in practice, and even if it is, the results~\eqref{eq:ff},\eqref{eq:flbad} still show SGN is stable, as the effect of roundoff errors are bounded independently of $\kappa_2(Y^T\!AX)$ and $\epsilon^{-1}$.} 

For the whole matrix, the computed version is such that the $i$th row is 
\[ [fl(AX(Y^T\!AX)_\epsilon^\dagger Y^TA)]_i = [AX]_i(Y^T\!AX+\epsilon_i)_\epsilon^\dagger Y^T\!A+\eepsilon
=[AX]_i(Y^T\!\At_iX)_\epsilon^\dagger Y^T\!A+\eepsilon
\]
for every $i$. In order to bound the associated error, we return to~\eqref{eq:sgnfirst} and note that the $i$th row of 
$A-[AX]_i(Y^T\!\At_iX)_\epsilon^\dagger Y^T\!A$ is that of $E_{\rm HMT}$ plus $AX(Y^T\!\At_i X)_\epsilon^\dagger Y^T U_\perp U_\perp^TAV_\perp V_\perp^T(I-\mathcal{P}_{X,V})=:AX(Y^T\!\At_i X)_\epsilon^\dagger E_2$. 
Denote by $W$ the $m\times n$ matrix whose $i$th row is $AX(Y^T\!\At_i X)_\epsilon^\dagger E_2$. 
Using the SVD $AX=U\Sigma V^T$, we have 
$[W]_i = [U]_i\Sigma V^T(Y^T\!\At_i X)_\epsilon^\dagger E_2$. 
Its norm can be bounded as
$\|[W]_i\|_2\leq \|[U]_i\|_2\|\Sigma V^T(Y^T\!\At_i X)_\epsilon^\dagger E_2 \|_2$. 
Using~\eqref{eq:specalways}, \eqref{eq:AXYnorm} 
and the assumption $\epsilon > \|Y^T\!(A-\At_i) X\|_2$, 
we obtain 
\begin{align*}
\|\Sigma V^T(Y^T\!\At_i X)_\epsilon^\dagger E_2 \|_2
\leq 2\|Y_1^\dagger\|_2\|E_2\|_2, 
\end{align*}
which holds for all $i$. 
Thus the Frobenius norm of $W$ is bounded by 
\begin{align}
\|W\|_F&  \leq \sqrt{\sum_{i}(\|[U]_i\|_2 2\|Y_1^\dagger\|_2\|E_2 \|_2)^2}   
= 2\sqrt{r}\|Y_1^\dagger\|_2\|E_2 \|_2. \label{eq:UYE}
\end{align}
We can now take the expectations with respect to $Y$ and $X$ separately as in~\eqref{eq:HMTF} to obtain 
  \begin{align*}
\mathbb{E}_{X,Y}\|W\|_F
&\leq 2\sqrt{r}
\sqrt{\mathbb{E}_Y\|Y_1^\dagger\|_2^2\ \mathbb{E}_X\|E_2 \|_2^2}
=2\sqrt{r}
\sqrt{\frac{4e^2(r+\ell)}{\ell^2}}
\sqrt{\mathbb{E}_X\|E_2 \|_2^2}\\
&\leq 2\sqrt{r}
\sqrt{\frac{4e^2(r+\ell)}{\ell^2}}
\sqrt{\mathbb{E}_X\|E_2 \|_F^2}
=2\sqrt{r}
\sqrt{\frac{4e^2(r+\ell)}{\ell^2}} \sqrt{r+\ell}\ \mathbb{E}_{X}\|E_{\rm HMT}\|_F\\
&=\frac{4e\sqrt{r}(r+\ell)}{\ell} \mathbb{E}\|E_{\rm HMT}\|_F.
  \end{align*}


\ignore{
\begin{align*}
\mathbb{E}\|E_2\|_2^2&=
\mathbb{E}\|Y^T U_\perp U_\perp^TAV_\perp V_\perp^T(I-\mathcal{P}_{X,V})\|_2^2\\
&\leq \mathbb{E}[\|Y^T U_\perp\|_2^2 \|U_\perp^TAV_\perp V_\perp^T(I-\mathcal{P}_{X,V})\|_2^2]\\
&= \mathbb{E}_{Y}\|Y_1^T\|^2\mathbb{E}_{X}\|M_X\|_2^2 
&\leq (\sqrt{r+\ell}+\sqrt{r})^2\mathbb{E}_{X}\|E_{\rm HMT}\|_2^2 \\
\end{align*}
so overall 
\[
\leq 
4e\sqrt{r}\frac{r+\ell}{\ell}\sqrt{\mathbb{E}_X[\|E_{\rm HMT}\|_2^2]}. 
\]
}

\ignore{
\begin{align}
\mathbb{E}&\|Y^T U_\perp U_\perp^TAV_\perp V_\perp^T(I-\mathcal{P}_{X,V})\|_2^2
=\mathbb{E}_{X,Y}\|Y^T M_X\|_2^2 
\leq (\sqrt{r+\ell}+\sqrt{r})^2\|M_X\|_F
(r+\ell)(\mathbb{E}_{X}\|M_X\|_F)^2 \nonumber \\
&= (r+\ell)(\mathbb{E}_{X}\|E_{\rm HMT}\|_F)^2 \leq \|\Sigma_2\|_F^2(r+\ell)(1+\frac{r}{r-\hat r-1}), \label{eq:HMTFG}
\end{align}
}


\qed

One might wonder what would happen if the assumption $\epsilon \geq \|Y^T\!(A-\At_i) X\|_2$ is violated; this is a condition that involves the unknowns $\At_i$. Fortunately, one can see that if we instead have 
$\epsilon \geq C\|Y^T\!(A-\At_i) X\|_2$ for some $C>0$, 
the bounds will largely remain the same; the bound on \eqref{eq:Efloat} would be multiplied by $C^{-1}$. Consequently, one can always safely take $\epsilon$ to be a small multiple of $u$, say $\epsilon=10u$. 

The bound~\eqref{eq:Efloat} is enough to show SGN is stable, as the effect of roundoff errors are bounded independently of $\kappa_2(Y^T\!AX)$ and $\epsilon^{-1}$, and $m,n$. However, \eqref{eq:Efloat} is larger than $E_{\rm GN}$ by a factor roughly $\sqrt{r}$, and this is an artifact of the analysis, namely the use of the loose inequality $\mathbb{E}[\|E_2 \|_2^2]\leq \mathbb{E}[\|E_2 \|_F^2]$ in the final inequality of the proof. 
To see that $\sqrt{r}$ should be removable, we note in~\eqref{eq:UYE} that $\|Y_1^\dagger\|_2\approx \frac{1}{\sqrt{r+\ell}-\sqrt{r}}$ by M-P with large deviation occuring with exponentially low probability, so $\mathbb{E}(\|Y_1^\dagger\|_2\|E_2 \|_2)\lesssim \frac{1}{\sqrt{r+\ell}-\sqrt{r}}\mathbb{E}\|E_2 \|_2$, 
and using~\cite[Prop.~A.2]{halko2011finding} we have 
$\mathbb{E}_{X}\|E_2\|_2\leq \sqrt{r+\ell}\|E_{\rm HMT}\|_2+\|E_{\rm HMT}\|_F$, indicating
\begin{equation} \label{eq:realfloat}  \mathbb{E}\|fl(AX(Y^T\!AX)_\epsilon^\dagger Y^TA) - A\|_F\lesssim 
\frac{\sqrt{r+\ell}+\sqrt{r}}{\sqrt{r+\ell}-\sqrt{r}}
\mathbb{E}\|E_{\rm HMT}\|_F
+\frac{2\sqrt{r(r+\ell)}}{\sqrt{r+\ell}-\sqrt{r}}\mathbb{E}\|E_{\rm HMT}\|_2.
\end{equation}

A key fact exploited in the proof of Theorem~\ref{thm:sta}, in particular~\eqref{eq:break},  is that for a linear system $Ax=b$, with a computed solution $\hat x$ satisfying $(A+\Delta A)\hat x=b$ we have $A\hat x=b+\eepsilon$ provided that $\|\hat x\|=\O(1)$, and the computd $x$ has much better accuracy than if the inverse $A^{-1}$ was computed explicitly. 
This phenomenon is mentioned in~\cite[Ch.14]{Higham:2002:ASNA}, and used in \cite{yamamotoetna2015}. 

\subsection{(In)stability of plain \Nystromm}\label{sec:insta}
It is natural to wonder, what could go wrong with the plain GN approximant~\eqref{eq:start} without the $\epsilon$-pseudoinverse? 
Two issues arise when one attempts to adapt the proof of Theorem~\ref{thm:sta}:
\begin{enumerate}
\item[(i)] The matrix $Y^T\!AX$ may not be numerically full rank, so \cite[Thm 21.4]{Higham:2002:ASNA} cannot be invoked. 
\item[(ii)] In the final step of the proof, the statement 
$\|[AX]_i(Y^T\!AX+\epsilon_i)^\dagger Y^T\!A\|=O(1)$ does not necessarily hold; recall the remark after Lemma~\ref{lem:AXX}.  Namely, without the $\epsilon$-pseudoinverse it is possible that 
$\|AX(Y^T\!AX+\epsilon_i)^\dagger\| \gg 1$, if 
$\|(Y^T\!AX+\epsilon_i)^\dagger\|$ happens to be $\gg 1/\epsilon$. 
This does not occur with high probability but one cannot rule out its possibility, as this depends on the behavior of roundoff errors. A further complication is that $\epsilon_i$ is not independent of $X,Y$, so bounds on $\mathbb{E}\|A X(Y^TA X+\epsilon_i)^\dagger\|_2^2$ as in~\eqref{eq:AX13} are not easy to obtain. 
\end{enumerate}

\ignore{
One can show that a computed version of GN results in $AX(Y^T\tilde AX)^{\dagger}Y^T\!A+O(\epsilon)$ (row-wise) for some $\tilde A$. 
The problem lies in~\eqref{eq:break}, where the fact 
$\|AX(Y^T\!AX+\epsilon_i)_\epsilon^\dagger\| =O(1)$ from~\eqref{eq:AX1} was used.} 

For a concrete example, suppose that $s_i$ as in 
\eqref{eq:mihatpert} was computed without the 
$\epsilon$-truncation so that the computed version is 
\[\tilde s_i = fl((Y^T\!AX+\epsilon_i)^\dagger [AX]_i^T),\]
implemented via the thin QR factorization 
$Y^T\!AX+\epsilon_i=QR$ as 
\[\tilde s_i = fl(R^{-1}(Q^T[AX]_i^T)).\]
As $Y^T\!AX$ is typically highly ill-conditioned, the bottom diagonal element of $R$ is $O(u)$; and there is nothing to stop it from being much smaller $\ll u$. If this happens, $\|\tilde s_i\|\gg 1$ and the argument in~\eqref{eq:break} breaks down, and one can see that this can give a poorly computed $fl(\hat A_r)$.
Otherwise, if $\|\tilde s_i\|=O(1)$, the issue (ii) is not present. 

Conversely, when neither (i) nor (ii) is present, the proof of Theorem~\ref{thm:sta} can be applied (and simplified) to show that plain GN is also stable. 
This suggests a natural and inexpensive way to implement a modified GN so that the outcome is a realization of an SGN, as we discuss in the next section. 

\section{Implementation}\label{sec:implement}
Here we discuss implementation details of (S)GN. The main topic is the pseudoinverse $(Y^T\!AX)^{\dagger}$, for which a careful implementation and modification make GN stable. We aim to find an efficient stabilization that avoids the $\epsilon$-truncation in the pseudoinverse unless necessary.
We also discuss a recommended choice of the oversampling parameter $\ell$. 
\subsection{Implementing the pseudoinverse}\label{sec:pinv}

The pseudoinverse $(Y^T\!AX)^\dagger$ is implemented via either 
(i) the QR factorization $Y^T\!AX=QR$, or (ii) the SVD. 
We mainly focus on the QR case as it is cheaper. 

The expression {\tt At = ((A*X)/(Y'*A*X))*(Y'*A)} does not always return stable results\footnote{This is because MATLAB's (back)slash commands {\tt /},{\tt $\backslash$} should be used with caution for underdetermined problems as they are not designed to find the minimum-norm solution via the QR factorization~\cite[\S 5.6]{Golubbookori}.
It is also important to note that the mathematically equivalent command {\tt At = ((A*X)*pinv(Y'*A*X))*(Y'*A)} is unstable and should not be used. 
}. 
An implementation of GN that works almost always 
  is to perform a QR factorization $Y^T\!AX=QR$ and take $((AX)R^{-1}) (Q^T(Y^T\!A))$, where $(AX)R^{-1}$ is obtained via triangular solve: 
in MATLAB, 
\begin{equation}  \label{eq:matlab}
{\tt AX = A*X;\quad YA = Y'*A;\quad [Q,R] = qr(Y'*AX,0);\quad At = (AX/R)*(Q'*YA)}  
\end{equation}

Note that the order in which the factors are multiplied is crucial for stability; for example, $(AX)(R^{-1}Q^T)(Y^T\!A)$ gives catastrophic results. 
Generally it is important that GN does not compute the core matrix $(Y^T\!AX)^{\dagger}$ explicitly, unlike other methods such as~\cite{tropp2017practical}. 

\subsubsection{SGN implementation}\label{sec:SGNimp}
To implement the stabilized GN so that the analysis in Section~\ref{sec:staGN} is valid, one approach is to compute the SVD of $Y^T\!AX$, truncate the singular values smaller than $\epsilon$, then apply the pseudoinverse. This reliable procedure costs $O(r^3)$ operations. 

A cheaper alternative is to compute the QR factorization $Y^T\!AX=QR$ as in \eqref{eq:matlab}, and look for diagonal element of $R$ less than $\epsilon$. Since $X,Y$ are random, this is a rank-revealing QR~\cite{chan1987rank,gu1996efficient} with high probability \cite{faststable}, and so by truncating the bottom-right corner of $R$, we obtain $Y^T\!AX=Q_1R_1+\eepsilon$ where $R_1$ is rectangular (fat) and upper triangular, and numerically full rank. One can further perform the QR of $R_1^T=Q_2R_2$ to form $Q_2R_2^{-1}Q_1^T$, which can be written as $(Y^T\!\At X)_\epsilon^\dagger$. We use this in our experiments. 

Another inexpensive hack is to perturb the computed $Y^T\!AX$ so that all singular values are larger than $u$; this way there is no need to truncate the singular values in the $\epsilon$-pseudoinverse. 
There are a few possible ways to achieve this: 
(a) compute the SVD and increase the small singular values, (b) compute a rank-revealing QR factorization, and perturb the diagonal elements of $R$ so they have $\geq u$ entries, and (c) simply perturb the diagonals of the (already computed) $R$. (a) and (b) are guaranteed to result in $\|R^{-1}\|=O(u^{-1})$. (c) is also expected to succeed with high probability. 
Each of these stabilization processes require only $O(r^3)$ operations, and works well in practice.

\subsubsection{Detecting potential instability in plain GN}\label{sec:fix}
The implementation~\eqref{eq:matlab} of plain GN is actually seen in experiments to be almost always stable. Indeed all experiments shown in this paper look almost identical between plain GN and stabilized GN, both in speed and error. It is therefore desirable to switch to SGN only when necessary. 

Recall the two issues discussed in Section~\ref{sec:insta}. 
Neither 
can happen unless $\|fl(R)^{-1}\|\gg u^{-1}$. This can be tested by a standard norm estimator for the computed $R$ based on the power method; each step costs only $O(r^2)$ operations. 
If one happens to have $\|fl(R)^{-1}\|\gg u^{-1}$, 
one switches to one of the implementations of SGN for an additional $O(r^3)$ cost, as described above. 


The above discussions also reveal that while plain GN can be unstable, instability is quite unlikely to actually manifest itself in practice; the situation is somewhat similar to LU with partial pivoting, which is unstable in the worst case but terrifically stable in practice~\cite[\S 22]{trefbau}, and is preferred to the QR-based linear solver, which is provenly backward stable but twice slower. 
Similarly, we suspect that in most cases, plain GN would be the preferred method for its simplicity and empirical stability, and our results explain why it should be stable most of the time. Perhaps this is the most practical message of this work: plain GN can be used without stability concerns most of the time, and if one wants to be sure, an inexpensive check and fix is available. 


\subsection{Convergence of GN and oversampling parameter $\ell$}\label{sec:oversampl}
What is a sensible choice of $\ell$? Our focus is to ensure the convergence of $\|A-A_r\|$ as $r$ grows, assuming that the singular values of $A$ decay, potentially quite slowly. 
A related discussion is given in \cite{tropp2017practical}, but the choice there requires the knowledge of the decay of $\sigma_i(A)$. 

The analysis of generalized \Nystrom 
indicates its optimality up to the factor $\sqrt{1+\frac{r}{r-\hat r-1}}\sqrt{1+\frac{r+\ell}{\ell-1}}$, or $\sqrt{1+\frac{r}{r-\hat r-1}}\frac{2\sqrt{e}(r+\ell)}{\ell}$ for SGN (which is likely a slight overestimate; one could also treat the estimate~\eqref{eq:realfloat} accounting for roundoff errors). 
Let us examine both terms in these products.

\begin{figure}[htpb]
  \begin{minipage}[t]{0.33\hsize}
      \includegraphics[height=43mm]{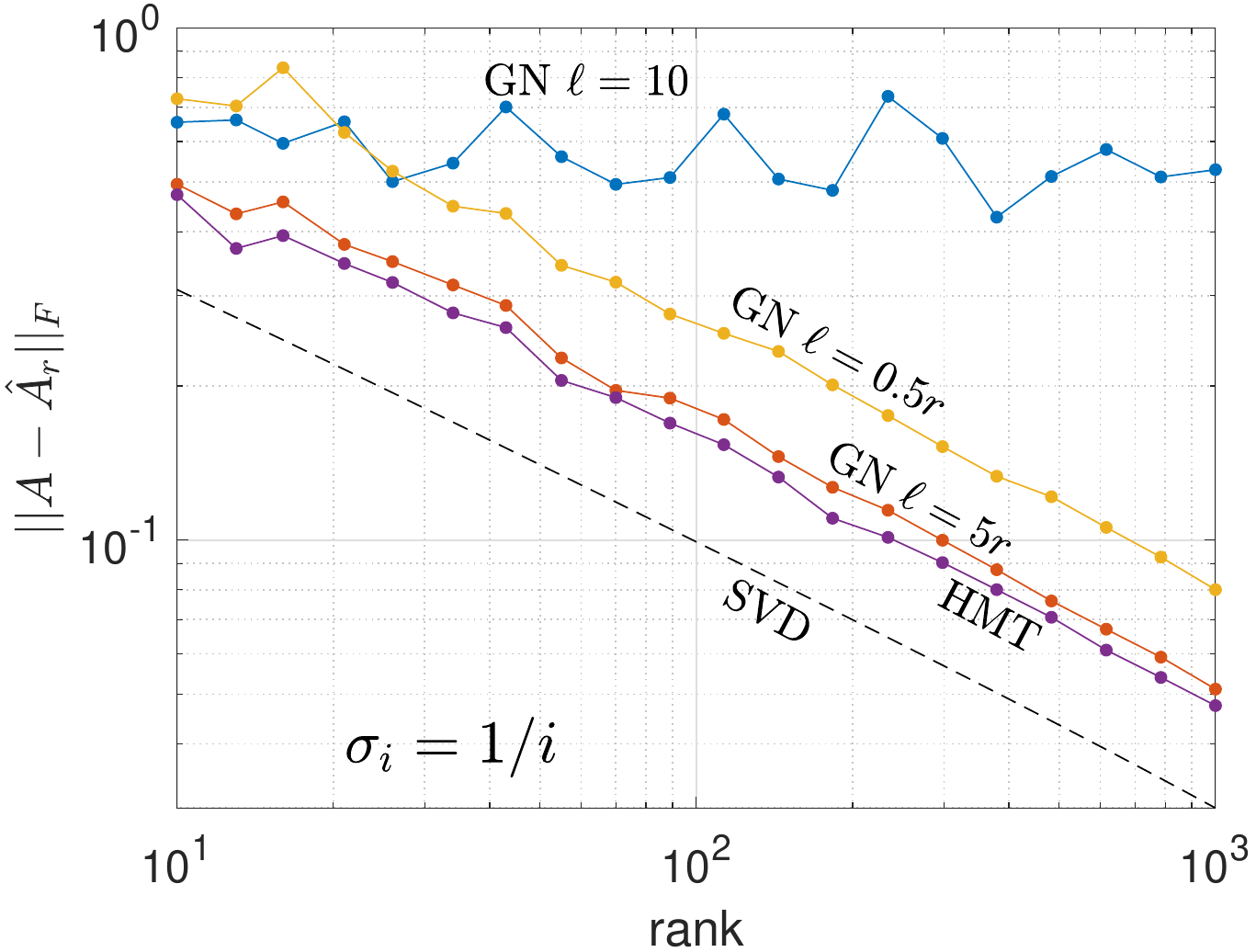}      

  \end{minipage}   
  \begin{minipage}[t]{0.33\hsize}
      \includegraphics[height=43mm]{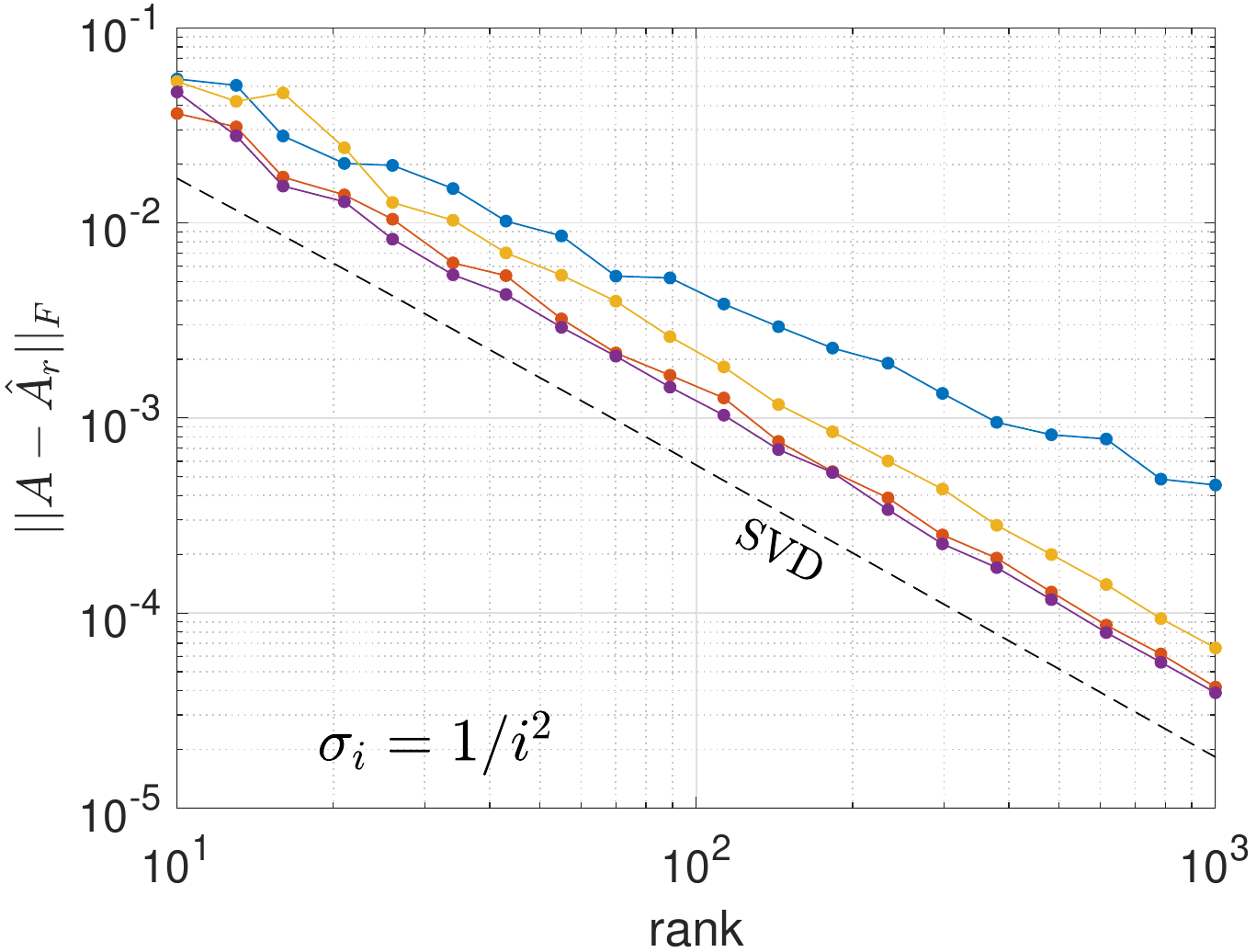}      
  \end{minipage}
  \begin{minipage}[t]{0.325\hsize}
      \includegraphics[height=43mm]{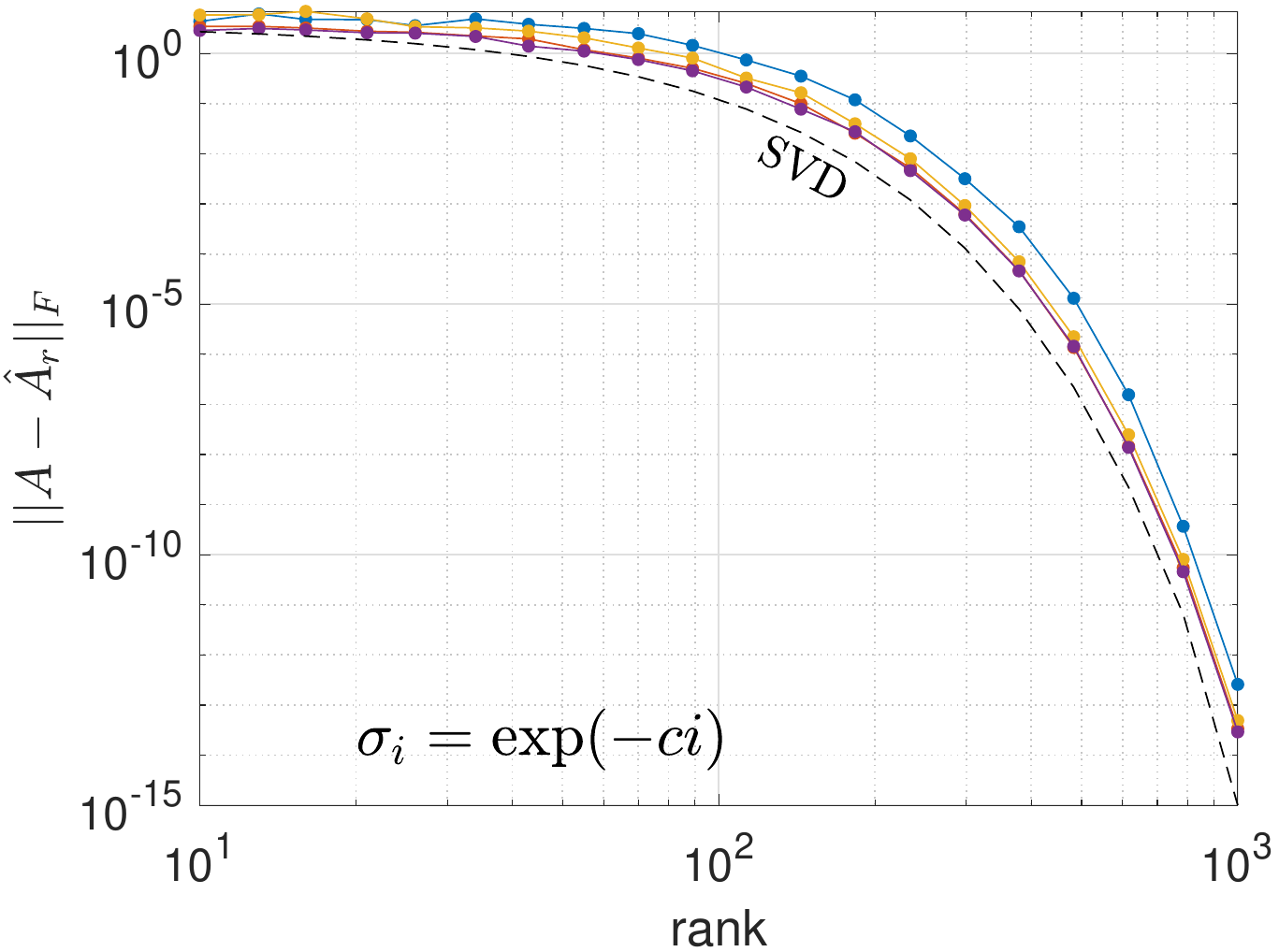}      
  \end{minipage}
  \caption{Generalized \Nystrom with different choices of oversampling parameter $\ell$: fixed $\ell=10$, recommended $\ell=r/2$, and large $\ell=5r$. 
Tested for matrices of varying spectrum: slow algebraic decay $\sigma_i(A)=1/i$ (left), algebraic decay $\sigma_i(A)=1/i^2$ (middle), 
and fast exponential decay 
$\sigma_i(A)=\exp(-ci)$ (right). 
}
  \label{fig:ell}
\end{figure}

The first term $\sqrt{1+\frac{r}{r-\hat r-1}}$ comes from the rangefinder and is independent of $\ell$. 
A practical consequence is that the approximation error of HMT (for which the second term is not present, recall~\eqref{eq:HMTbound}) always traces that of the optimal truncated SVD. 
It is worth noting that while many references suggest a constant oversampling suffices, e.g. $r=\hat r+5$ (which is fine if $r=O(1)$), for $\sqrt{1+\frac{r}{r-\hat r-1}}$ to be constant we need $r-\hat r$ to scale linearly with $r$. This is illustrated in Figure~\ref{fig:ell} (right), where the horizontal gap between SVD and HMT is roughly constant in log-scale. 


The second term $\sqrt{1+\frac{r+\ell}{\ell-1}}$ or $\frac{2\sqrt{e}(r+\ell)}{\ell}$  highlights the importance of the choice of $\ell$. 
For example, suppose we fix $\ell$ to be constant, say 10, as we grow $r$. Then these terms grow with $r$. In fact, if the decay of $\sigma_i(A)$ is not very fast, $\|E_{\rm GN}\|$ could even grow as we increase $r$, which is clearly undesirable. 
This is illustrated in Figure~\ref{fig:ell}, where fixing $\ell$ is evidently not enough when $\sigma_i(A)=1/i$. The issue improves when $\sigma_i(A)$ decay faster (e.g. exponential decay or $\sigma_i(A)=1/i^s$ with $s\geq 2$), but the fact remains that GN with fixed $\ell$ gets farther from optimal as $r$ increases, whereas choosing $\ell=cr$ avoids this issue. 

Based on this discussion, to ensure the approximation error of GN traces that of truncated SVD, a reliable choice is to let $\ell$ grow proportionally to $r$. Then
these terms 
 are bounded by a constant independent of $r$. 
Figure~\ref{fig:ell} shows how taking $\ell=cr$ results in GN converging proportionally to $\sigma_i(A)$; larger $c$ gets GN closer to HMT, but this obviously comes with more computation and storage. In all other experiments we use our recommended choice $\ell=r/2$.


\section{Resampling and updating the matrix}\label{sec:update}
The simplicity of GN makes it amenable and flexible in a variety of situations. One example is resampling: if one finds the approximant $\hat A_r$ to be insufficient in accuracy, a standard remedy is to increase the rank to $r+\delta r$; the associated cost is $O(N_{\delta r})$. Alternative methods would require orthogonalization, so GN is much more economical. 

For another example, suppose that data has been appended, so that one has 
the matrix $\begin{bmatrix}  A\\ B\end{bmatrix}$ with $B\in\mathbb{R}^{\hat m\times n}$, where we have a low-rank approximation to $A\approx AX(Y^T\!AX)_\epsilon^{\dagger}Y^T\!A $. 
Finding a low-rank approximation to $\begin{bmatrix}  A\\ B\end{bmatrix}$ is a simple matter of computing $BX$ and $\tilde Y^TB$ where $\tilde Y$ is a new sketch matrix, and then 
\[
\begin{bmatrix}  A\\ B\end{bmatrix}
\approx 
\begin{bmatrix}
AX\\BX  
\end{bmatrix}
\left([Y^T, \tilde Y^T]
\begin{bmatrix}  A\\ B\end{bmatrix}
X\right)_\epsilon^{\dagger}[Y^T\!A+\tilde Y^TB]
=\begin{bmatrix}
AX\\BX  
\end{bmatrix}
\left(
Y^T\!AX+ \tilde Y^TBX\right)_\epsilon^{\dagger}[Y^T\!A+\tilde Y^TB]. 
\]
The computational cost in updating the approximation is $O(\hat mn\log n+r^3)$, assuming $B$ is dense. Clearly, appending columns $[A,B]$ can be handled analogously. 

This is a notable advantage of GN---other algorithms cannot deal with such updates nearly as efficiently, because the orthogonalization step will have to be recomputed, involving $O((m+n)r^2)$ operations. 

Another example of update that we can handle (which other sketching/streaming algorithms can also do), 
 is when the matrix undergoes perturbation $A\leftarrow A+E$, or the streaming model. In this case we can update the sketches $AX,Y^T\!A$ accordingly via $AX\leftarrow AX+EX, Y^T\!A\leftarrow Y^T\!A+Y^TE$. Often $E$ has structure (e.g. sparsity) that allows us to compute $EX,Y^TE$ efficiently. 

\section{Experiments}\label{sec:ex}
We report numerical experiments to illustrate the performance of GN in comparison with other algorithms. We set $\epsilon=10^{-15}$ in SGN, and  
in all methods we take $X,Y$ to be subsampled DCT matrices. 
\subsection{Dense matrices}
We first address large dense matrices. We take $A=U\Sigma V^T\in\mathbb{R}^{50000\times 50000}$ with geometrically decaying singular values, where $U,V$ are obtained via the QR factorization of square Gaussian matrices. 
For $r=1000,2000,\ldots, 10^4$, we compute rank-$r$ approximations $A_k\approx A$ via HMT~\cite{halko2011finding}, Tropp~\cite{tropp2017practical} and GN, and compare the runtime and the approximation quality $\|A-A_k\|_F$. 
This is the same setup as in Figure~\ref{fig:posdef}, but for a nonsymmetric matrix. 

\begin{figure}[htpb]
  \begin{minipage}[t]{0.5\hsize}
      \includegraphics[height=50mm]{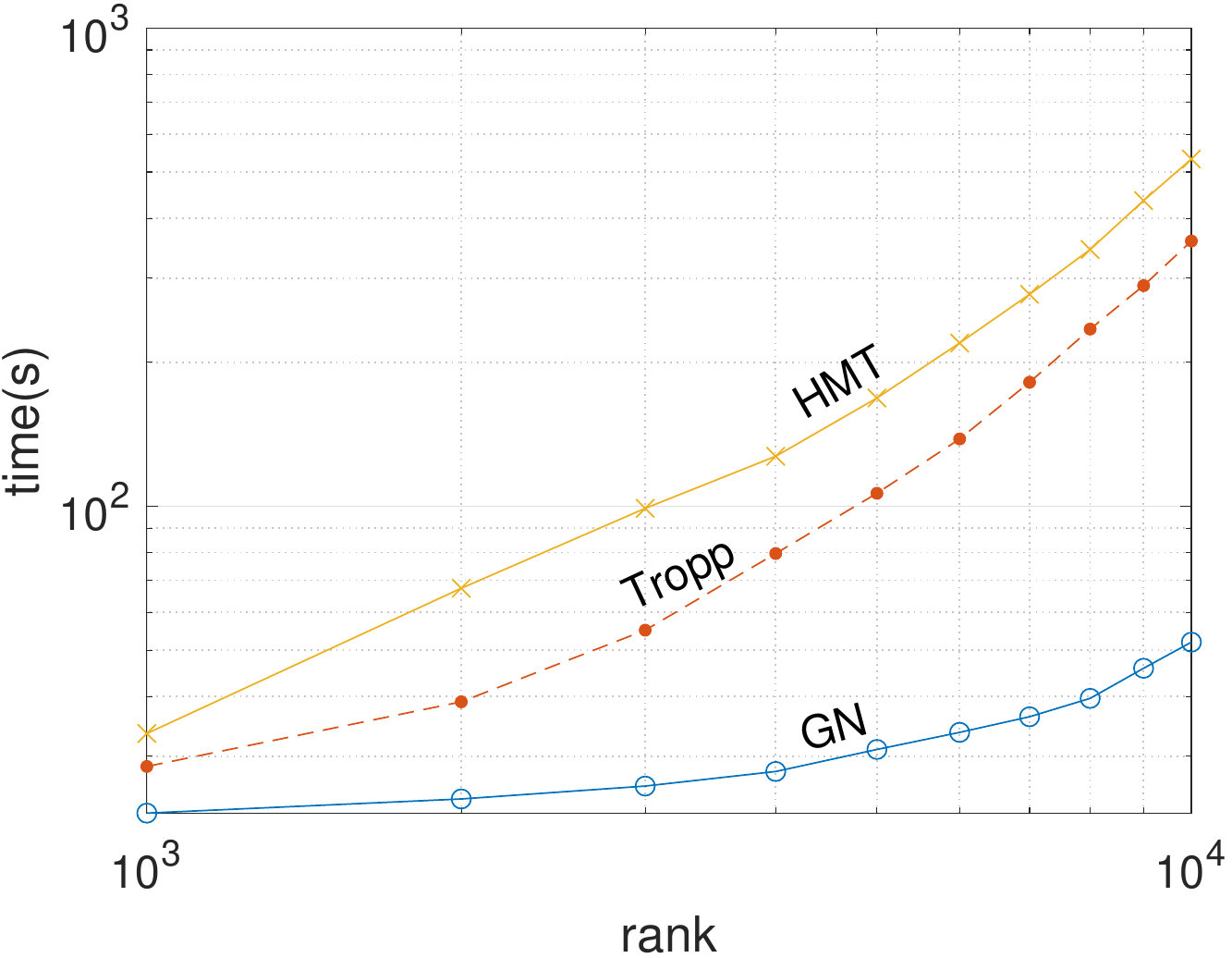}      
  \end{minipage}   
  \begin{minipage}[t]{0.5\hsize}
      \includegraphics[height=50mm]{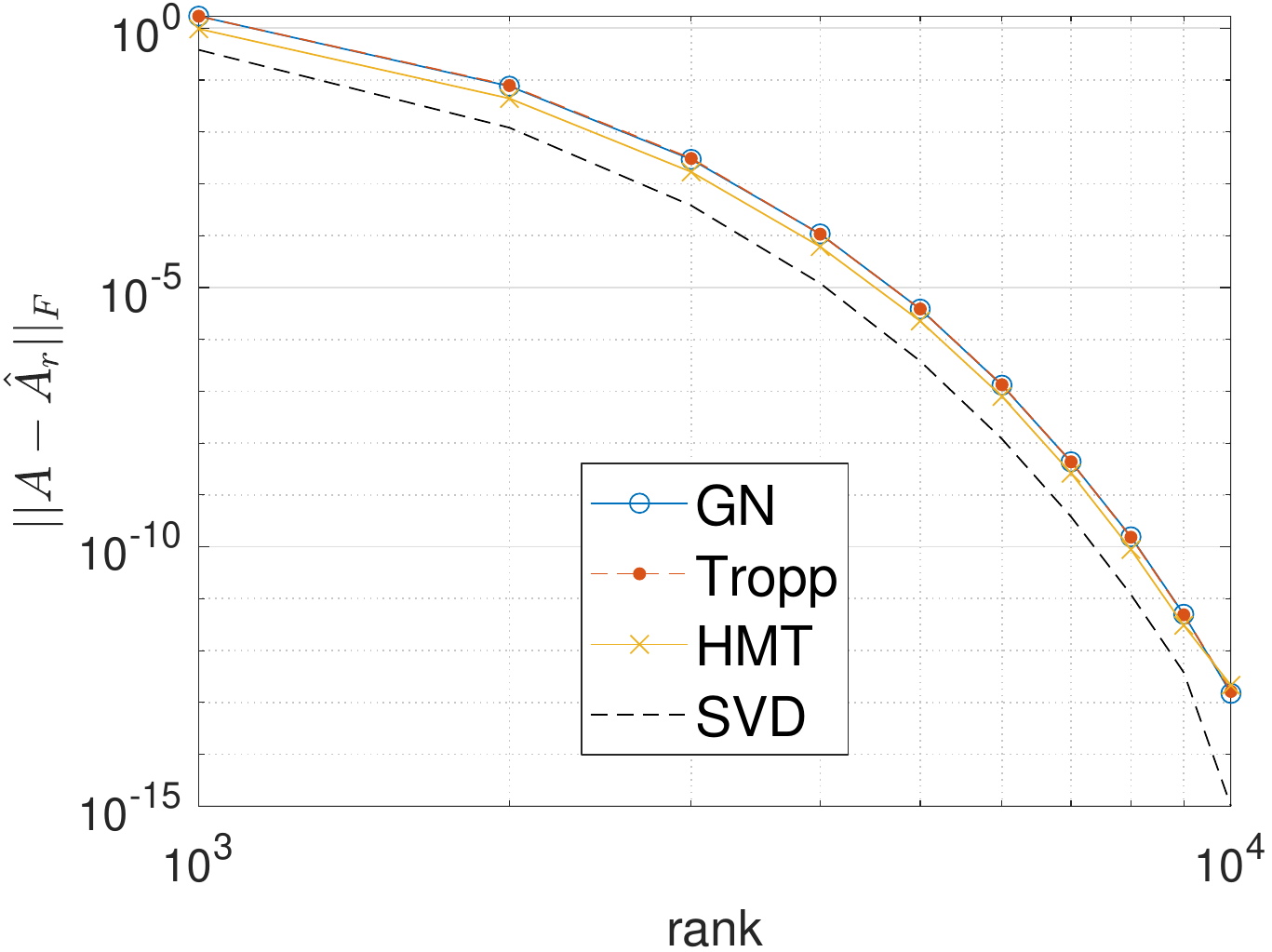}            
  \end{minipage}
  \caption{Comparison of algorithms on a $50000\times 50000$ dense matrix with geometrically decaying singular values.}
  \label{fig:compalgs}
\end{figure}

We observe that GN is significantly faster than HMT and Tropp, the speedup increasing with $r$, as expected. 

This experiment also illustrates the numerical stability established in Section~\ref{sec:gennyst}: Recall that Tropp and GN are mathematically equivalent when roundoff errors are ignored. The fact that the accuracy of the two methods are nearly identical in Figure~\ref{fig:compalgs} verifies that roundoff errors do not negatively affect GN. 

The equivalence between GN and Tropp up to numerical errors, together with the stability analysis that shows numerical errors have negligible effects on GN,
implies that the observations made in the extensive experiments reported in \cite{tropp2017practical} in terms of the excellent accuracy (but not the speed) of Tropp apply also to GN. 
\subsection{When should a randomized algorithm be used?}
A natural question is: 
when should we choose randomized algorithms over classical, deterministic algorithms? It is clear from the complexity that when $r\ll m,n$, randomized algorithms would outperform classical methods that need $O(mn^2)$ cost. However, when $r\approx O(m,n)$, it is unclear if randomized algorithms are still competitive. 

\begin{figure}[htpb]
  \begin{minipage}[t]{0.5\hsize}
      \includegraphics[height=50mm]{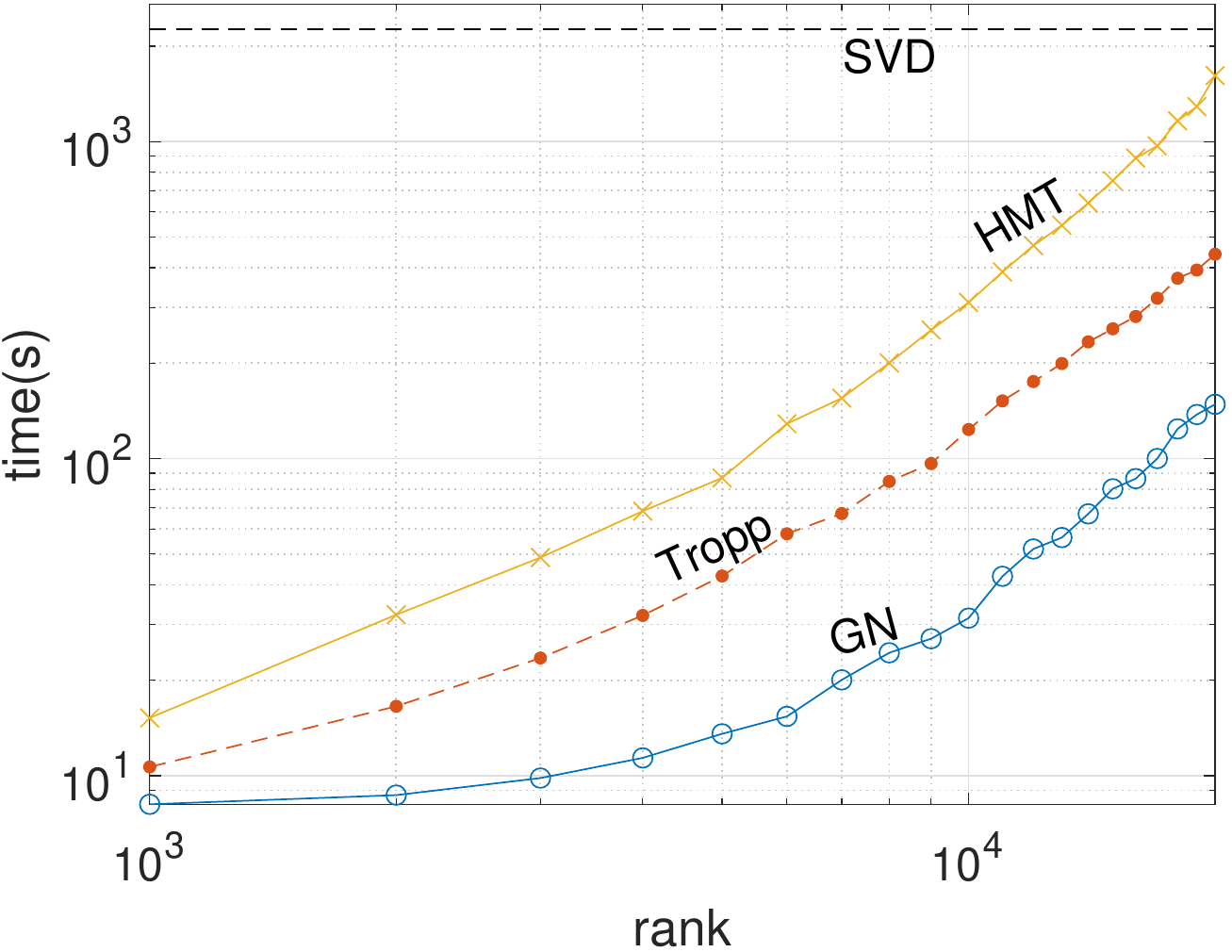}
  \end{minipage}   
  \begin{minipage}[t]{0.5\hsize}
      \includegraphics[height=50mm]{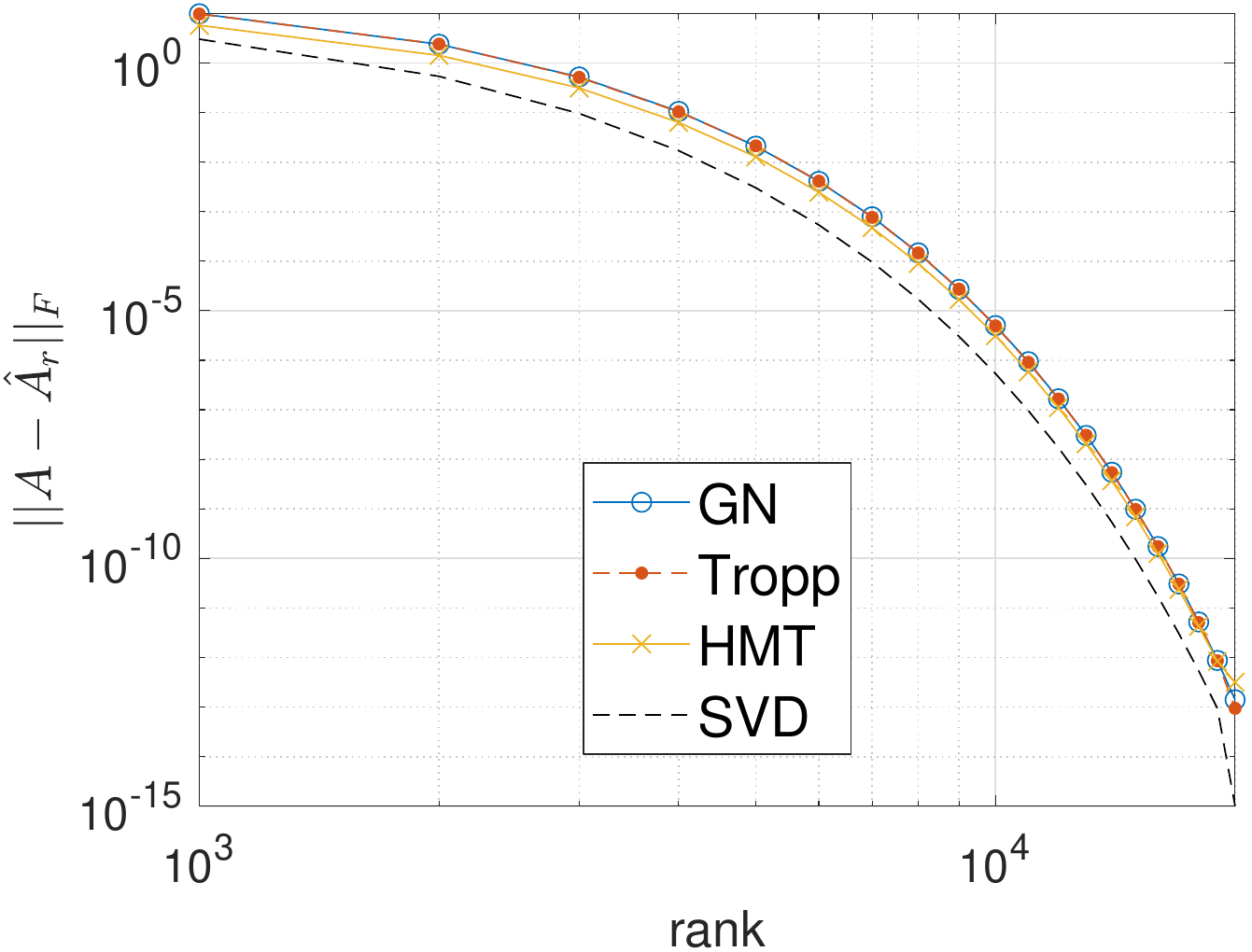}
  \end{minipage}
  \caption{Performance compared with a full SVD.
The matrix is $30000\times 30000$ dense, with geometrically decaying singular values.
 }
  \label{fig:vsfull}
\end{figure}

To gain insight, here we compare randomized algorithms with MATLAB's SVD, which performs bidiagonalization followed by divide-and-conquer. We form a random $30000\times 30000$ matrix with geometrically decaying singular values (but with slow decay; $\sigma_1=1,\sigma_{20000}=10^{-15}$). 
We vary the required rank $r$ from $1000$ to $20000=\frac{2}{3}n$. 
The runtime and approximation accuracy are shown in Figure~\ref{fig:vsfull}. 

We see that randomized algorithms outperform classical algorithms in efficiency by a large margin, even when the rank is quite close to $m,n$. In particular, a 10-fold speedup is observed with GN even when $r=\frac{2}{3}n$; note that at this point, there is no saving in memory over storing the entire matrix $A$, as $Y^TA$ is square with the recommended oversampling factor $1.5$. 
For this reason $r=\frac{2}{3}n$ is a practical 'limit' with GN. 

In terms of accuracy, classical SVD of course finds the exact truncated SVD (up to an $O(u)$ backward error). Randomized algorithms are optimal to within small factors; the difference is usually offset by taking a slightly larger $r$ (as in Figure~\ref{fig:vsfull}), but more oversampling would be required when the singular values decay slowly, as illustrated in Figure~\ref{fig:ell}. 

The main message here is that randomized algorithms appear to be preferable whenever the matrix admits an storage-efficient low-rank approximation. 

\subsection{Sparse matrices}\label{sec:exsparse}
For sparse matrices the cost $N_r$ in Table~\ref{tab:tab} reduces to $N_r\approx \mbox{nnz}(A)r$, and the relative importance of the orthogonalization cost $O(mr^2)$ can increase. 
We take the matrix "transient" from~\cite{davis2011university}, which is of size $m=n=178,866$ and $961,368$ nonzero entries. We again see significant speedup with GN.

\begin{figure}[htpb]
  \begin{minipage}[t]{0.5\hsize}
      \includegraphics[height=50mm]{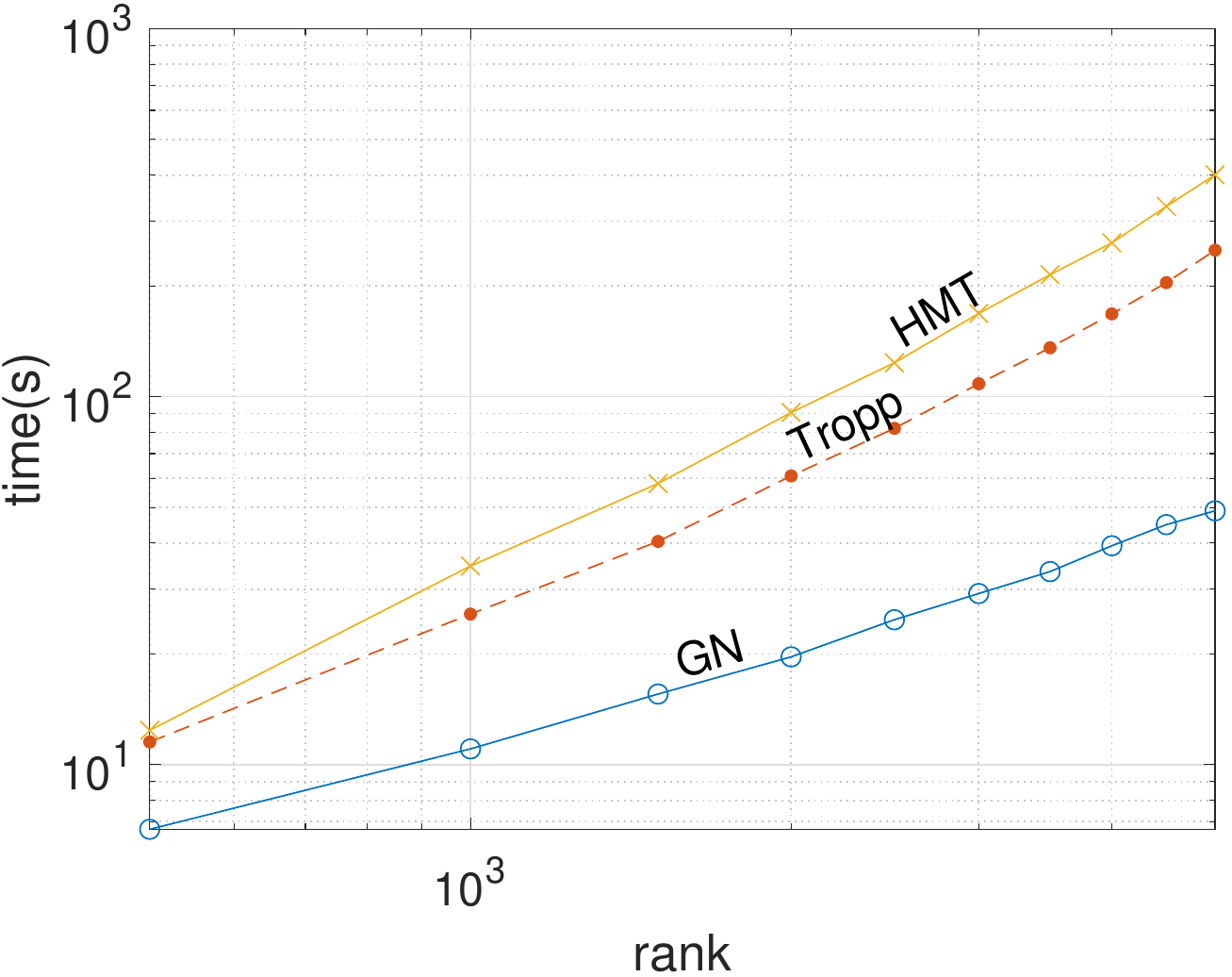}
  \end{minipage}   
  \begin{minipage}[t]{0.5\hsize}
      \includegraphics[height=50mm]{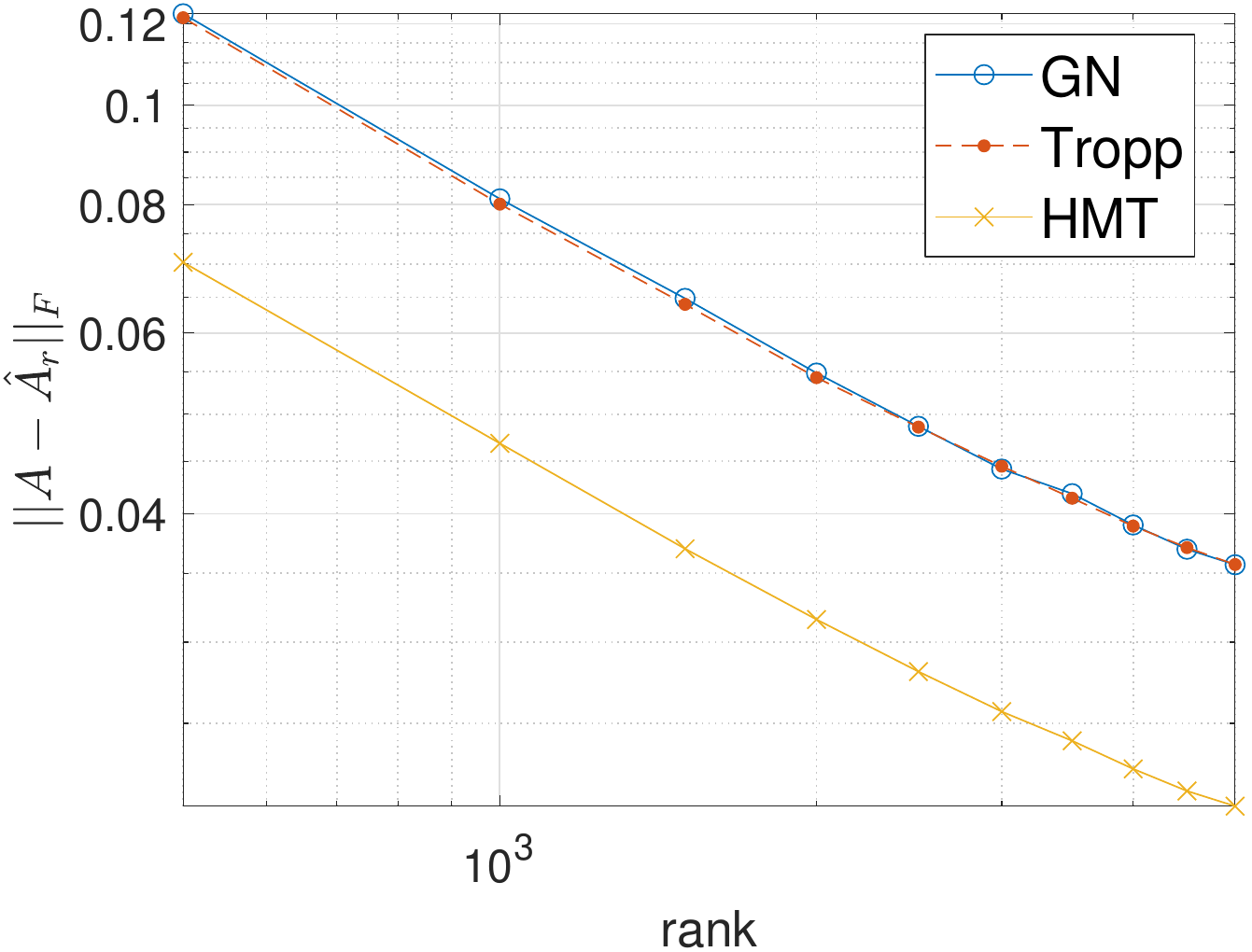}
  \end{minipage}
  \caption{Performance for a sparse matrix.}
  \label{fig:sparse}
\end{figure}

As mentioned in Section~\ref{sec:cost}, for sparse matrices, whether the low-rank approximation here is effective relative to storing the nonzero elements is debatable. Indeed it is sensible to force the approximants to also have sparse factors, as e.g. in the CUR factorization or~\cite{grigori2018low,zhang2002low}; 
our GN does not have this property. One could take $X,Y$ to be sparse sketch matrices, such as the Countsketch matrix; this then reduces to the algorithm by Clarkson and Woodruff~\cite{clarkson2009numerical}. 

\subsection{Matrix updating}
\begin{figure}[htpb]
  \begin{minipage}[t]{0.5\hsize}
      \includegraphics[height=50mm]{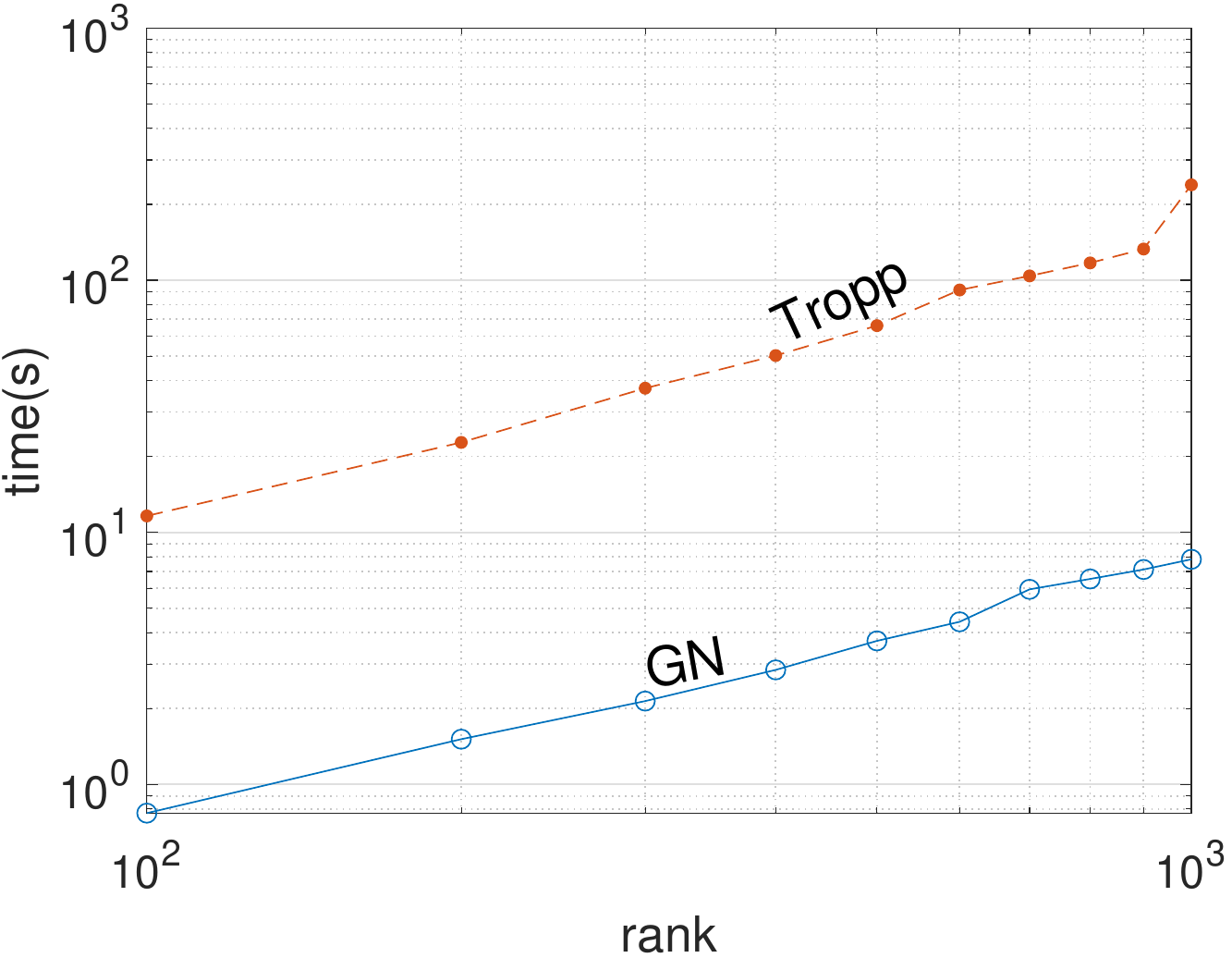}      
  \end{minipage}   
  \begin{minipage}[t]{0.5\hsize}
      \includegraphics[height=50mm]{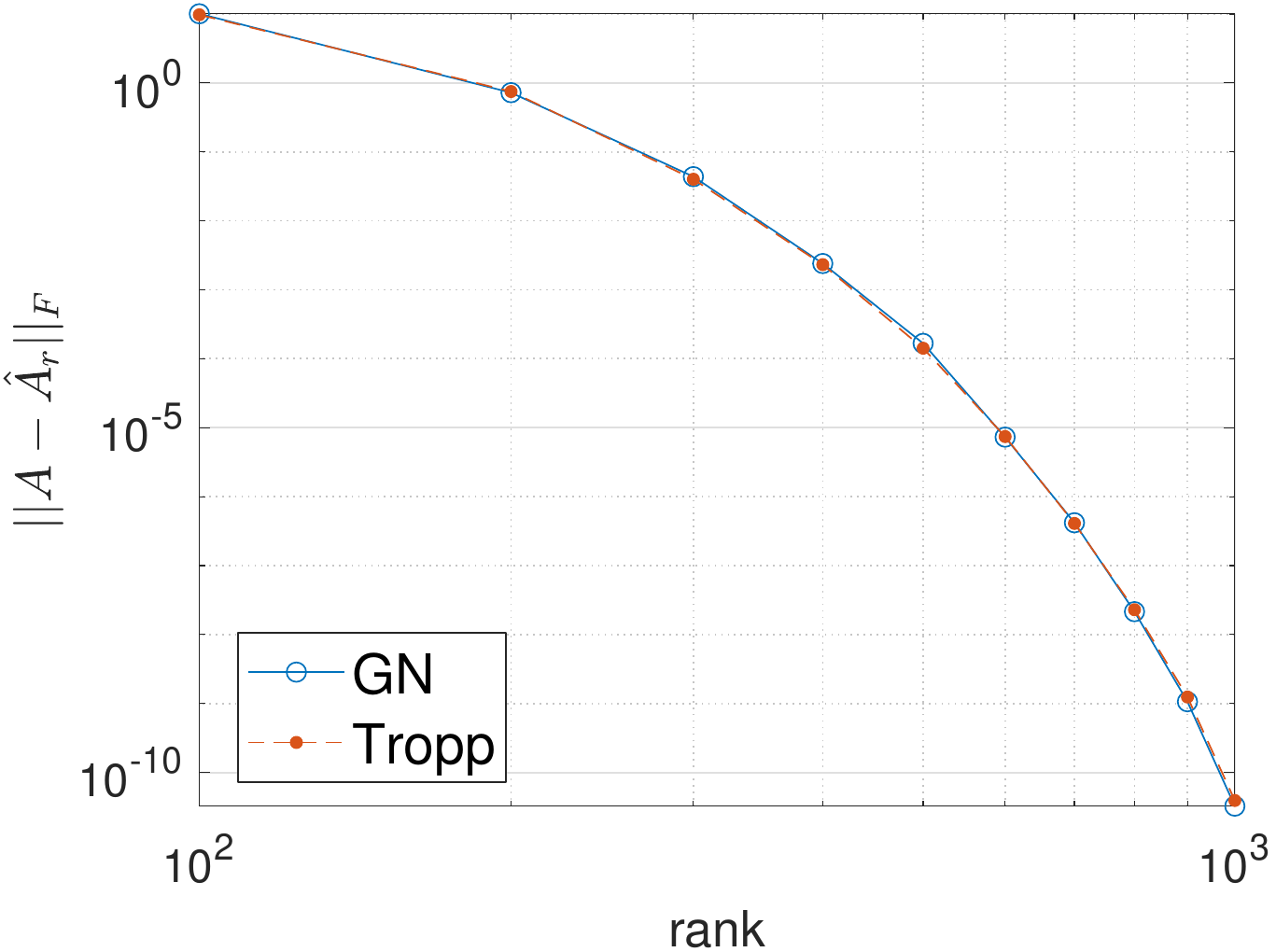}            
  \end{minipage}
  \caption{Updating the matrix. 
}
  \label{fig:append}
\end{figure}

Finally, 
we start with a $2\cdot 10^6\times 2\cdot 10^6$ sparse matrix with geometrically decaying singular values generated by $A=A_1DA_2^T$, where $A_1,A_2$ are sparse and $D$ is diagonal with geometrically decaying entries. We then append 1000 rows at the bottom $A\leftarrow \big[
\begin{smallmatrix}
  A\\ E
\end{smallmatrix}
\big]$. The update in GN is performed as described in Section~\ref{sec:update}. 
Here we only compare GN and Tropp, as HMT would require revisiting the original matrix $A$, which is not suitable in this situation. The results are in Figure~\ref{fig:append}. Once again, GN achieves significant speedup by avoiding orthogonalization. 

\subsection*{Acknowledgment}
I would like to thank Gunnar Martinsson and Maike Meier for helpful discussions and pointers to the literature, and Jared Tanner for suggesting the experiment in Figure~\ref{fig:vsfull}.




\def\noopsort#1{}\def\l{\char32l}\def\v#1{{\accent20 #1}}
  \let\^^_=\v\def\hbk{hardback}\def\pbk{paperback}

\end{document}